\documentclass[10pt]{article}
\usepackage{amsmath, amssymb,epsfig,bm}
\usepackage{latexsym}

\usepackage[hypertex,hyperindex]{hyperref}
\usepackage{graphicx}
\usepackage{color}
\numberwithin{equation}{section}
\DeclareMathAlphabet{\itbf}{OML}{cmm}{b}{it}

\newcommand{\PP}{\mathbb{P}}
\newcommand{\EE}{\mathbb{E}}
\newcommand{\RR}{\mathbb{R}}

\newcommand{\ds}{\displaystyle}

\def\debproof{ {\bf Proof.} }
\def\finproof{\hfill {\small $\Box$} \\}

\newcommand{\bp}{{\itbf p}}
\newcommand{\bq}{{\itbf q}}
\newcommand{\br}{{\itbf r}}
\newcommand{\bs}{{\itbf s}}
\newcommand{\bz}{{\itbf z}}
\newcommand{\bx}{{\itbf x}}
\newcommand{\bv}{{\itbf v}}
\newcommand{\bev}{{\itbf e}}
\newcommand{\bw}{{\itbf w}}
\newcommand{\bu}{{\itbf u}}
\newcommand{\bn}{{\itbf n}}
\newcommand{\by}{{\itbf y}}
\newcommand{\bX}{{\bf X}}

\newcommand{\bR}{{\itbf R}}
\renewcommand{\i}{\mathrm{i}}

\newcommand{\bE}{{\itbf E}}
\newcommand{\bH}{{\itbf H}}
\newcommand{\bJ}{{\itbf J}}
\newcommand{\bF}{{\itbf F}}

\newcommand{\be}{\begin{eqnarray}}
\newcommand{\ee}{\end{eqnarray}}
\newcommand{\nn}{\nonumber}
\newcommand{\ben}{\begin{eqnarray*}}
\newcommand{\een}{\end{eqnarray*}}

\def\ds{\displaystyle}
\def\nm{\noalign{\medskip}}

\newtheorem{lem}{Lemma}[section]
\newtheorem{prop}{Proposition}[section]
\newtheorem{cor}{Corollary}[section]
\newtheorem{thm}{Theorem}[section]

\newtheorem{definition}{Definition}[section]

\newcommand{\rot}{{\nabla\times}}

\newcommand{\ddiv}{{\nabla\cdot}}

\begin{document}

\title{Target Detection and Characterization from Electromagnetic
Induction Data\thanks{\footnotesize This work was supported  by
ERC Advanced Grant Project MULTIMOD--267184, China NSF under the
grants 11001150, 11171040, and 11021101, and National Basic
Research Project under the grant 2011CB309700.}}

\author{Habib Ammari\thanks{\footnotesize Department of Mathematics and Applications,
Ecole Normale Sup\'erieure, 45 Rue d'Ulm, 75005 Paris, France
(habib.ammari@ens.fr).} \and Junqing Chen\thanks{\footnotesize
Department of Mathematical Sciences, Tsinghua University, Beijing
100084, China (jqchen@math.tsinghua.edu.cn).} \and Zhiming
Chen\thanks{\footnotesize LSEC, Institute of Computational
Mathematics, Chinese Academy of Sciences, Beijing 100190, China
(zmchen@lsec.cc.ac.cn).} \and
 \and Josselin Garnier\thanks{\footnotesize
Laboratoire de Probabilit\'es et Mod\`eles Al\'eatoires \&
Laboratoire Jacques-Louis Lions, Universit\'e Paris VII, 75205
Paris Cedex 13, France (garnier@math.univ-paris-diderot.fr).}  \and Darko
Volkov\thanks{Department of Mathematical Sciences, Stratton Hall,
100 Institute Road, Worcester, MA 01609-2280, USA
(darko@wpi.edu).}}

\maketitle

\begin{abstract}
The goal of this paper is to contribute to the field of
nondestructive testing by eddy currents. We provide a mathematical
analysis and a numerical framework for simulating the imaging of
arbitrarily shaped small volume conductive inclusions from
electromagnetic induction data. We derive, with proof, a
small-volume expansion of the eddy current data measured away from
the conductive inclusion. The formula involves two polarization
tensors: one associated with the magnetic contrast and the second
with the conductivity of the inclusion.  Based on this new
formula, we design a location search algorithm. We include in this
paper a discussion on data sampling, noise reduction, and on
 probability of detection. We provide numerical examples that
 support our findings.
\end{abstract}

\bigskip

\noindent {\footnotesize Mathematics Subject Classification
(MSC2000): 35R30, 35B30}


\section{Introduction}

Nondestructive testing by eddy currents is a  technology of choice
in the assessment of  the structural integrity of a variety of
materials such as, for instance, aircrafts or  metal beams, see
\cite{review}. It is also of interest in technologies related to
safety of public arenas where a large number of people have to be
screened.

We introduce in this paper a novel analysis pertaining to
small-volume expansions for eddy current equations, which we then
apply to developing new imaging techniques. Our mathematical
analysis extends recently established results and methods for full
Maxwell's equations to the eddy current regime.

We propose a new eddy current reconstruction method relying on the
assumption that the objects to be imaged are small. This present
study is related to the theory of small-volume perturbations of
Maxwell's equations, see \cite{AVV}. It is, however, specific to
eddy currents  and to the particular lengthscales relevant to that
case.

We first note that in the eddy current regime a diffusion equation is used for modeling
electromagnetic fields. The characteristic length is the skin
depth of the conductive object to be imaged \cite{review}.
%
We consider the regime where the skin depth is comparable to the
characteristic size  of the conductive inclusion.


Using the $\bE$-formulation for  the eddy current problem, we
first establish energy estimates. We start from integral
representation formulas for the electromagnetic fields arising in
the presence of a small conductive inclusion
 to
 rigorously derive an asymptotic
expansion for   the  magnetic part of the  field. The effect of
the conductive targets on the magnetic field measured away from
the target is expressed  in terms of two polarization tensors: one
associated with the magnetic contrast (called magnetic
polarization tensor) and the second with the conductivity of the
target (called conductivity polarization tensor). The magnetic
polarization tensor has been first introduced in \cite{AVV} in the
zero conductivity case while the concept of conductivity
polarization tensor appears to be new.

Based on our asymptotic formula we are then able to construct a
localization method for  the conductive inclusion. That method
involves  a response matrix data. A MUSIC (which stands for
MUltiple Signal Classification) imaging functional is proposed for
locating the target. It uses the projection of a magnetic dipole
located at the search point onto the image space of the response
matrix. Once the location is found, geometric features of the
inclusion may be reconstructed using a least-squares method. These
geometric features together with material parameters (electric
conductivity and magnetic permeability) are incorporated in the
conductivity and magnetic polarization tensors. It is worth
emphasizing that, as will be shown by our asymptotic expansion,
the perturbations in the magnetic field due to the presence of the
inclusion are complex-valued while the unperturbed field can be
chosen to be real. As consequence, we only process the imaginary
part of the recorded perturbations. Doing so, we do not need to
know the unperturbed field with an order of accuracy higher than
the order of magnitude of the perturbation. An approximation of
lower order of the unperturbed field is enough.

The so called Hadamard measurement sampling technique is applied
in order to reduce the impact of noise in measurements. We briefly
explain some underlying basic ideas. Moreover, we provide
statistical distributions for the singular values of the response
matrix in the presence of measurement noise. An important strength
of our analysis is that it can be applied for rectangular response
matrices. Finally, we simulate our localization technique on a
test example.

The paper is organized as follows. Section 2 is devoted to a
variational formulation of the eddy current equations. Section 3
contains the main contributions of this paper. It provides a
rigorous derivation of the effect of a small conductive inclusion
on the magnetic field measured away from the inclusion. Section 4
extends MUSIC-type localization proposed in \cite{AILP} to the
eddy current model. Section 5 discusses the effect of noise on the
inclusion detection and proposes a detection test based on the
significant eigenvalues of the response matrix. Section 6
illustrates numerically on test examples our main findings in this
paper. A few concluding remarks are given in the last section.

\section{Eddy Current Equations}

Suppose that there is an electromagnetic inclusion in $\RR^3$ of
the form $B_\alpha=\bz+\alpha B$, where $B \subset \RR^3$ is a bounded, smooth
 domain containing the origin. Let $\Gamma$ and $\Gamma_\alpha$ denote the boundary of $B$
  and $B_\alpha$. Let $\mu_0$ denote the magnetic permeability
of  the free space. Let $\mu_*$ and $\sigma_*$
denote the permeability  and the conductivity of the  inclusion which are also assumed to be
constant. We introduce the piecewise constant magnetic
permeability and electric conductivity
\begin{eqnarray*}
\mu_{\alpha}(\bx) = \left \{
    \begin{array}{ll}
    \mu_*  & \mbox  {in $B_\alpha$,   } \\
\nm
        \mu_{0}   & \mbox  {in $B^c_\alpha=\RR^3 \backslash  \bar{B}_\alpha,$}
           \end{array}
 \right.\ \ \ \
 \sigma_\alpha (\bx)=  \left \{
    \begin{array}{ll}
    \sigma_*  & \mbox  {in $B_\alpha$,  }\\
\nm
        0   & \mbox  {in $B_\alpha^c$.}
           \end{array}
 \right.
\end{eqnarray*}

Let $(\bE_\alpha, \bH_\alpha)$ denote the eddy current fields in
the presence of the electromagnetic inclusion $B_\alpha$ and a
source current $\bJ_0$ located outside the inclusion. Moreover, we
suppose that $\bJ_0$ has a compact support and is divergence free:
$\nabla \cdot \bJ_0 =0$ in $\RR^3$. The fields $\bE_\alpha$ and
$\bH_\alpha$ are the solutions of the following eddy current
equations: \be \label{alpha} \left\{
\begin{array}{l}
\rot \bE_\alpha =\i\omega \mu_\alpha \bH_\alpha \ \ \
\hbox{ in } \RR^3,
\\
\nm
 \rot  \bH_\alpha =~\sigma_\alpha \bE_\alpha + \bJ_0 \ \ \
 \hbox{ in } \RR^3,\\
 \nm
 \bE_\alpha(\bx)=O(|\bx|^{-1}), \ \ \bH_\alpha(\bx)=O(|\bx|^{-1}) \ \ \ \mbox{ as } |\bx|\rightarrow \infty.
\end{array} \right.
\ee By eliminating $\bH_\alpha$ in \eqref{alpha} we obtain the
following $\bE$-formulation of the eddy current problem
\eqref{alpha}: \be \label{eform} \left\{
\begin{array}{l}
\rot\mu^{-1}_\alpha\rot\bE_\alpha-
\i\omega\sigma_\alpha\bE_\alpha= \i\omega\bJ_0\ \ \  \mbox{ in }
\RR^3, \\ \nm \ddiv\bE_\alpha =0\ \ \ \mbox{ in } B_\alpha^c, \\
\nm \bE_\alpha(\bx)=O(|\bx|^{-1})\ \ \ \mbox{ as }
|\bx|\rightarrow\infty.
\end{array}\right.
\ee Throughout this paper, let $\bu^\pm$ denote the limit values
of $\bu(\bx \pm t \bn)$ as $ t \searrow 0$, where $\bn$ is the outward normal to
$\Gamma_\alpha$, if they exist.
 We will use the function spaces
$$\bX_\alpha(\RR^3)=\left\{\bu:\frac{\bu}{\sqrt{1+|\bx|^2}}\in L^2(\RR^3)^3,
\rot\bu\in L^2(\RR^3)^3, \ddiv\bu=0 \mbox{ in } B_\alpha^c
\right\},
$$
and
$$\widetilde{\bX}_\alpha(\RR^3)=\left\{\bu: {\bu} \in \bX_\alpha(\RR^3),
\int_{\Gamma_\alpha} \bu^+\cdot \bn =0  \right\},
$$
and the sesquilinear form on $\widetilde{\bX}_\alpha(\RR^3)\times
\widetilde{\bX}_\alpha(\RR^3)$
$$a_\alpha(\bE, \bv)=(\mu^{-1}_\alpha\rot\bE,\rot\bv)_{\RR^3}- \i\omega \sigma_*(\bE,\bv)_{B_\alpha},$$
where $(\cdot,\cdot)_D$ stands for the $L^2$ inner product on the
domain $D\subseteq \RR^3$. The weak formulation of the
$\bE$-formulation (\ref{eform}) is: Find $\bE_\alpha\in
\widetilde{\bX}_\alpha(\RR^3)$ such that \be \label{weformtilde}
a_\alpha(\bE_\alpha,\bv)= \i\omega(\bJ_0,\bv)_{B^c_\alpha},\ \ \
\forall\bv\in \widetilde{\bX}_\alpha(\RR^3). \ee The uniqueness
and existence of solution of the problem \eqref{weformtilde} is
known (cf., e.g., Ammari et al. \cite{ABN} and Hiptmair
\cite{hiptmair}). Note that the constraint $\int_{\Gamma_\alpha}
\bu^+\cdot \bn =0$ in $\widetilde{\bX}_\alpha(\RR^3)$ only serves
to enforce the uniqueness of $\bE_\alpha$ in $B^c_\alpha$
\cite{hiptmair}. This is not essential for the validity of the
$\bE$-formulation of the eddy current model. We have \be
\label{weform} a_\alpha(\bE_\alpha,\bv)=
\i\omega(\bJ_0,\bv)_{B^c_\alpha},\ \ \ \forall\bv\in
{\bX}_\alpha(\RR^3). \ee

Throughout the paper we denote by  $\bE_0$ the unique solution of the problem
\begin{equation} \label{E0}
\left\{
\begin{array}{l}
\rot\mu^{-1}_0\rot\bE_0=\i\omega\bJ_0\ \ \  \mbox{ in } \RR^3 , \\
\nm \ddiv \bE_0 = 0 \mbox{ in } \RR^3 , \ \ \  \\ \nm
\bE_0(\bx)=O(|\bx|^{-1})\ \ \  \mbox{ as } |\bx| \rightarrow
\infty.
\end{array}\right.
\end{equation}
The field $\bE_0$  satisfies \be
\label{e0}
(\mu^{-1}_0\rot\bE_0,\rot\bv)_{\RR^3}=\i\omega(\bJ_0,\bv)_{\RR^3},
\ \ \ \forall\bv\in {\bf H}_{-1}(\mathbf{curl};\RR^3), \ee where
$\displaystyle{\bf
H}_{-1}(\mathbf{curl};\RR^3)=\Big\{\bu:\frac{\bu}{\sqrt{1+|\bx|^2}}\in
L^2(\RR^3)^3, \rot\bu\in L^2(\RR^3)^3\Big\}$.

\section{Derivation of the Asymptotic Formulas}

In this section we will derive the asymptotic formula for
$\bH_\alpha$ when the inclusion is small. Let
$k=\omega\mu_0\sigma_*$. We are interested in the asymptotic regime
when $\alpha \rightarrow 0$ and
\begin{equation} \label{deff5}
\nu = k \alpha^2
\end{equation}
is of order one. Moreover, we assume that $\mu_*$ and $\mu_0$ are
of the same order.

In eddy current testing the wave equation is converted into the
diffusion equation, where the characteristic length is the skin
depth $\delta$, given by $\delta = \sqrt{2/k}$. Hence, in the
regime $\nu =O(1)$, the skin depth $\delta$ is of order of the
characteristic size $\alpha$ of the inclusion.

We will always denote by $C$ a
 generic constant which depends possibly on $\mu_*/\mu_0$, the upper bound of
$\omega \mu_0\sigma_* \alpha^2$,
 the domain $B$, but is independent of
 $\omega,\sigma_*,\mu_0,\mu_*$. Let $\mu_r  =\mu_*/\mu_0.$

\subsection{Energy Estimates}

We start with the following estimate.

\begin{lem}\label{L:3.1}
There exists a constant $C$ such that
$$\|\rot(\bE_\alpha-\bE_0)\|_{L^2(\RR^3)}+
\sqrt k\|\bE_\alpha-\bE_0\|_{L^2(B_\alpha)}\leq C\alpha^{3/2} (
\sqrt k
\|\bE_0\|_{L^\infty(B_\alpha)}+\|\rot\bE_0\|_{L^\infty(B_\alpha)}).$$
\end{lem}

\debproof
By \eqref{weform} and \eqref{e0}, we know that
\begin{eqnarray}\label{x1}
& &(\mu^{-1}_\alpha\rot(\bE_\alpha-\bE_0),\rot\bv)_{\RR^3}- \i\omega(\sigma_\alpha(\bE_\alpha-\bE_0),\bv)_{B_\alpha}\nn\\
\nm &
&\qquad=(\mu_0^{-1}-\mu_*^{-1})(\rot\bE_0,\rot\bv)_{B_\alpha}+
\i\omega(\sigma_\alpha\bE_0,\bv)_{B_\alpha},\ \ \ \forall
\bv\in\bX_\alpha(\RR^3).
\end{eqnarray}
Since
$$|(\rot\bE_0,\rot\bv)_{B_\alpha}|\leq C\alpha^{3/2}\|\rot\bE_0\|_{L^\infty(B_\alpha)}\|\rot\bv\|_{L^2(B_\alpha)}$$
and
$$|(\sigma_\alpha\bE_0,\bv)|\leq C\alpha^{3/2}\sigma_*\|\bE_0\|_{L^\infty(B_\alpha)}\|\bv\|_{L^2(B_\alpha)},$$
by taking $\bv = \bE_\alpha-\bE_0\in\bX_\alpha(\RR^3)$ in
\eqref{x1} and multiplying the obtained equation by $\mu_0$ we
have that \ben
& &\mu_r^{-1}\| \rot (\bE_\alpha-\bE_0)
\|_{L^2(\RR^3)}^2+k\|\bE_\alpha-\bE_0\|_{L^2(B_\alpha)}^2\\
\nm
&\le&C\alpha^{3/2}\left(\|\rot\bE_0\|_{L^\infty(B_\alpha)}\|\rot(\bE_\alpha-\bE_0)\|_{L^2(B_\alpha)}+k\|\bE_0\|_{L^\infty(B_\alpha)}\|\bE_\alpha-\bE_0\|_{L^2(B_\alpha)}\right).
\een This completes the proof. \finproof
\bigskip

Let $H^1(B_\alpha)=\{ \varphi \in L^2(B_\alpha), \nabla \varphi
\in L^2(B_\alpha)^3\}$.
Let $\phi_0\in H^1(B_\alpha)$ be the
solution of the problem
\be
\label{phi0} -\Delta\phi_0 =-\nabla\cdot\bF \mbox{ in
} B_\alpha,\ \ \ -\partial_\bn\phi_0 = (\bE_0
(\bx)-\bF(\bx))\cdot\bn \mbox{ on } \Gamma_\alpha,\ \ \
\int_{B_\alpha}\phi_0 \,d\bx=0,
\ee
where
\begin{equation}
\label{deff}
\bF(\bx)= \frac12
[\nabla_\bz \times\bE_0(\bz) ]\times(\bx-\bz)+\frac 13 [ {\bf D}_\bz ( \nabla_\bz \times \bE_0) (\bz)](\bx-\bz)\times(\bx-\bz).
\end{equation}
Here $[{\bf D}_\bz  (\nabla_\bz \times \bE_0)(\bz)]_{ij} = \partial_{z_i} (\nabla_\bz \times \bE_0(\bz))_j$
 is the $(i,j)$-th element of the gradient matrix ${\bf D}_\bz  (\nabla_\bz \times \bE_0)(\bz)$ of $\nabla_\bz\times\bE_0(\bz)$.
Let ${\rm tr}$ denote the trace. Since ${\rm
tr}[{\bf D} (\nabla  \times \bE_0)]=\nabla\cdot(\nabla\times\bE_0)=0$, we
know that
\begin{equation}
 \label{deff2}
\nabla\times\bF(\bx)=\nabla_\bz\times\bE_0(\bz)+[{\bf D}_\bz(\nabla_\bz\times\bE_0)(\bz)](\bx-\bz).
\end{equation}
Note that since $\bE_0$ is smooth in $\bar{B}_\alpha$ we have
\begin{equation} \label{deff3}
\|\nabla\times \bE_0 - \nabla\times \bF\|_{L^\infty(B_\alpha)}\le
C\alpha^2\|\nabla\times\bE_0\|_{W^{2,\infty}(B_\alpha)}.
\end{equation}

Denote by $\displaystyle \bH_0= (\i\omega\mu_0)^{-1}\rot\bE_0$ and
introduce $\bw\in \widetilde{\bX}_\alpha(\RR^3)$ as the solution
of the problem \be \label{correct_w} a_\alpha(\bw,\bv)& =
&\i\omega\mu_0(\mu_0^{-1}-\mu_*^{-1})(\bH_0(\bz) + {\bf D}
\bH_0(\bz) (\bx-\bz),\rot\bv)_{B_\alpha} \nonumber \\ \nm && +
\i\omega (\sigma_\alpha\bF,\bv)_{B_\alpha}, \qquad \forall\bv\in
\widetilde{\bX}_\alpha. \ee
The following lemma provides a higher-order correction of the
error estimate in Lemma \ref{L:3.1}.

\begin{lem}\label{L:3.2}
Let $\bw$ be defined by \eqref{correct_w}. There exists a constant
$C$ such that
 \begin{eqnarray}
 & &\|\rot(\bE_\alpha-\bE_0-\bw)\|_{L^2(\RR^3)}  \leq C \alpha^{7/2}(|1-\mu_r^{-1}|+ \nu)\|\nabla\times\bE_0\|_{W^{2,\infty}(B_\alpha)},\label{x2} \\
 & &\|\bE_\alpha-\bE_0-\nabla\phi_0-\bw\|_{L^2(B_\alpha)}
\leq C \alpha^{9/2}(|1-\mu_r^{-1}|+
\nu)\|\nabla\times\bE_0\|_{W^{2,\infty}(B_\alpha)},\label{x3}
 \end{eqnarray}
 where $\nu$ is given by \eqref{deff5}.
\end{lem}

\debproof First we set $\psi \in H^1(B_\alpha)$ and $g$ be such
that $g=\psi$ on $\Gamma_\alpha$, $\Delta g= 0$ in $B_\alpha^c$,
and $g = O(|\bx|^{-1})$ at infinity. Let $$\bv :=
\left\{\begin{array}{l} \nabla \psi \quad \mbox{in } B_\alpha,\\
\nm \nabla g \quad \mbox{in } B_\alpha^c. \end{array} \right.$$
Since $\bv \in \bX_\alpha$, it follows from (\ref{weform}) that
$$
\i
\omega(\sigma_\alpha\bE_\alpha, \nabla\psi)_{B_\alpha}=0,\ \ \ \forall \psi\in H^1(B_\alpha).
$$
This yields $\ddiv\bE_\alpha=0$ in $B_\alpha$ and $\bE_\alpha^-
\cdot\bn=0$ on $\Gamma_\alpha$.

Similarly, we can show from (\ref{correct_w}) that $\bw^-\cdot\bn=
- \bF(\bx) \cdot \bn$ on $\Gamma_\alpha$ and
$\ddiv\bw=-\nabla\cdot\bF$ in $B_\alpha$. From \eqref{phi0} we
also know that $\ddiv(\bE_0+\nabla\phi_0)=\nabla\cdot\bF$ in
$B_\alpha$ and $(\bE_0+\nabla\phi_0)^-\cdot\bn= \bF(\bx)\cdot\bn$
on $\Gamma_\alpha$. Thus \ben
\nabla\cdot\left(\bE_\alpha-\bE_0-\nabla\phi_0-\bw\right)=0\ \
\mbox{in }B_\alpha,\ \ \ \
\left(\bE_\alpha-\bE_0-\nabla\phi_0-\bw\right)^-\cdot\bn=0\ \
\mbox{on }\Gamma_\alpha, \een which implies by scaling argument
and the embedding theorem that \ben
\|\bE_\alpha-\bE_0-\nabla\phi_0-\bw\|_{L^2(B_\alpha)}&\le&C\alpha\|\nabla\times\left(\bE_\alpha-\bE_0-\nabla\phi_0-\bw\right)\|_{L^2(B_\alpha)}\\
&=&C\alpha\|\nabla\times\left(\bE_\alpha-\bE_0-\bw\right)\|_{L^2(B_\alpha)},
\een for some constant $C$ independent of $\alpha$ and $\sigma_*$.
Therefore, \eqref{x3} follows from \eqref{x2}.

To show \eqref{x2}, we define $\tilde{\phi}_0$ as the solution of
the exterior problem
\begin{equation*}
-\Delta\tilde{\phi}_0=0 \mbox{ in } B^c_\alpha,\ \ \  \tilde{\phi}_0 = \phi_0 \mbox{ on } \Gamma_\alpha,\ \ \  \tilde{\phi}_0\rightarrow 0 \mbox{ as }|\bx|\rightarrow\infty.
\end{equation*}
The existence of $\tilde{\phi}_0$ in $\displaystyle
W^{1,-1}(B^c_\alpha)=\Big\{\varphi:\frac{\varphi}{\sqrt{1+|\bx|^2}}\in
L^2(B_\alpha^c), \nabla \varphi \in L^2(B_\alpha^c)^3 \Big\}$ is
known (cf., e.g., N\'ed\'elec \cite{nedelec}).

Define ${\bm \Phi}_0=\nabla\phi_0$ in $B_\alpha$, ${\bm
\Phi}_0=\nabla\tilde\phi_0$ in $B_\alpha^c$, then ${\bm \Phi}_0\in
\bX_\alpha(\RR^3)$. It follows from \eqref{x1} and
\eqref{correct_w} that for all $\bv \in \bX_\alpha(\RR^3)$
\begin{eqnarray*}
(\mu^{-1}_\alpha\rot(\bE_\alpha-\bE_0-{\bm
\Phi}_0-\bw),\rot\bv)_{\RR^3}- \i\omega(\sigma_\alpha
(\bE_\alpha-\bE_0-{\bm \Phi}_0-\bw),\bv)_{B_\alpha} \nonumber\\
\nm =\i\omega\mu_0(\mu_0^{-1}-\mu_*^{-1})(\bH_0-\bH_0(\bz) - {\bf
D} \bH_0(\bz) (\bx-\bz),\rot\bv)_{B_\alpha}
+\i\omega(\sigma_\alpha (\bE_0+{\bm \Phi}_0 -\bF),\bv)_{B_\alpha}.
\end{eqnarray*}
By multiplying the above equation by $\mu_0$ we have then
\be\label{correct_all}
(\mu_0\mu_\alpha^{-1}\rot(\bE_\alpha-\bE_0-{\bm
\Phi}_0-\bw),\rot\bv)_{\RR^3}-\i k(\bE_\alpha-\bE_0-{\bm
\Phi}_0-\bw,\bv)
_{B_\alpha}\nonumber\\
=\i\omega\mu_0(1-\mu_r^{-1})(\bH_0 -\bH_0(\bz)- {\bf D} \bH_0(\bz)
(\bx-\bz),\rot\bv)_{B_\alpha}+\i k(\bE_0+{\bm \Phi}_0
-\bF,\bv)_{B_\alpha}. \ee It is easy to check that \ben
|\i\omega\mu_0(\bH_0-\bH_0(\bz) - {\bf
D} \bH_0(\bz) (\bx-\bz),\rot\bv)_{B_\alpha}|&\leq& C\alpha^{7/2}\|\i\omega\mu_0\bH_0\|_{W^{2,\infty}(B_\alpha)}\|\rot\bv\|_{L^2(B_\alpha)}\\
&=&C\alpha^{7/2}\|\nabla\times\bE_0\|_{W^{2,\infty}(B_\alpha)}\|\rot\bv\|_{L^2(B_\alpha)}.
\een Now taking $\bv=\bE_\alpha-\bE_0-{\bm
\Phi}_0-\bw\in\bX_\alpha(\RR^3)$ in \eqref{correct_all}, since
$\nabla\times{\bm \Phi}_0=0$ in $\RR^3$ and ${\bm
\Phi}_0=\nabla\phi_0$ in $B_\alpha$, we obtain that
\begin{eqnarray*}
&&\|\rot(\bE_\alpha-\bE_0-\bw)\|^2_{L^2(\RR^3)}+k\|\bE_\alpha-\bE_0-\nabla\phi_0-\bw\|^2_{L^2(B_\alpha)}\\
\nm &\leq&
C\alpha^{7/2}|1-\mu_r^{-1}|\|\nabla\times\bE_0\|_{W^{2,\infty}(B_\alpha)}\|\rot\bv\|_{L^2(B_\alpha)}
+ k \|\bE_0 - \bF +\nabla\phi_0 \|_{L^2(B_\alpha)}\|\bv\|_{L^2(B_\alpha)}\\
\nm &\leq& C \alpha^{7/2} (|1-\mu_r^{-1}|+ \nu)
\|\nabla\times\bE_0\|_{W^{2,\infty}(B_\alpha)} \|\rot \bv\|_{L^2(
B_\alpha)}.
\end{eqnarray*}
Here, we have used
\begin{equation}\label{deff4} \|\bE_0 -\bF +\nabla\phi_0\|_{L^2(B_\alpha)}\le
C\alpha\|\nabla\times(\bE_0-\bF)\|_{L^2(B_\alpha)}\le
C\alpha^{9/2}\|\nabla\times\bE_0\|_{W^{2,\infty}(B_\alpha)}
\end{equation}
and $\|\bv\|_{L^2(B_\alpha)}\le
C\alpha\|\nabla\times\bv\|_{L^2(B_\alpha)}$, since $\bE_0 - \bF
+\nabla\phi_0$ and $\bv$ are divergence free in $B_\alpha$ and
have vanishing normal traces on $\Gamma_\alpha$. This shows
\eqref{x2} and completes the proof.
\finproof

We note that ${\bf D} \bH_0(\bz)$ is symmetric since $\nabla
\times \bH_0(\bz)=0$. Hence, by Green's formula, \ben &&\ds
(\mu_0^{-1}-\mu_*^{-1})(\bH_0(\bz) + {\bf D} \bH_0(\bz)
(\bx-\bz),\rot\bv)_{B_\alpha}\\
\nm &=&\ds (\mu_0^{-1}-\mu_*^{-1})\int_{\Gamma_\alpha}(
(\bH_0(\bz) + {\bf D} \bH_0(\bz)
(\bx-\bz))\times\bn)\cdot\bv d\bx \\
\nm &=&\ds \int_{\Gamma_\alpha}[\mu_\alpha^{-1}(\bH_0(\bz) + {\bf
D} \bH_0(\bz) (\bx-\bz))\times\bn]_{\Gamma_\alpha}\cdot\bv d \bx,
\een where $[\cdot]_{\Gamma_\alpha}$ stands for the jump of the
function across $\Gamma$. Let $\hat\bw({\bm
\xi})=\bw(\bz+\alpha{\bm \xi})$, we know from \eqref{correct_w}
that, $\forall\bv\in \widetilde{\bX}_1(\RR^3)$,
\begin{eqnarray*}
\ds
(\mu^{-1}\rot\hat\bw,\rot\bv)_{\RR^3}-\i\omega\alpha^2(\sigma\hat\bw,\bv)_B
&=&\ds \i \alpha \omega\mu_0 \int_{\Gamma}[\mu^{-1}(\bH_0(\bz) +
\alpha {\bf D} \bH_0(\bz) \bm \xi)\times\bn]_{\Gamma}\cdot\bv
d{\bm \xi}
\\ && \ds +\i \omega\alpha^2 (\sigma \bF(\bz+\alpha\bm\xi),\bv)_B ,
\end{eqnarray*}
where $\mu({\bm \xi})=\mu_*$ if ${\bm \xi}\in B$, $\mu({\bm \xi})=\mu_0$ if ${\bm \xi}\in B^c$ and
$\sigma({\bm \xi})=\sigma_*$ if ${\bm \xi}\in B$, $\sigma({\bm \xi})=0$ if ${\bm \xi}\in B^c$.

This motivates us to introduce the solution $\bw_0({\bm \xi})$ of the interface problem
\begin{equation}
\label{defwo} \left\{
\begin{array}{l}
\nabla_{\bm \xi} \times \mu^{-1}  \nabla_{\bm \xi} \times \bw_0- \i\omega\sigma\alpha^2\bw_0= \i \omega
\sigma \alpha^2 \left[\alpha^{-1} \bF(\bz+\alpha\bm\xi)\right] \mbox{ in } B\cup B^c, \\
\nm \nabla_{\bm \xi} \cdot \bw_0 = 0 \mbox{ in } B^c, \\ \nm
\mbox{[}\bw_0\times\bn ]_\Gamma= 0, \ \ \ [\mu^{-1} \nabla_{\bm
\xi} \times \bw_0\times\bn]_\Gamma= -\i\omega (1 - \mu_r^{-1})
(\bH_0(\bz) +
\alpha {\bf D} \bH_0(\bz) \bm \xi)\times\bn \mbox{ on } \Gamma,\\
\nm \bw_0(\bm \xi)= O(|\bm\xi|^{-1}) \mbox{ as } |\bm
\xi|\rightarrow\infty.
\end{array}\right.
\end{equation} It is easy to check that
$\displaystyle\bw(\bx)=\alpha\bw_0\Big(\frac{\bx-\bz}{\alpha}\Big)$.

The following theorem which is the main result of this section now
follows directly from Lemma \ref{L:3.2}.

\begin{thm}\label{T:3.1}
There exists a constant $C$ such that \ben &
&\Big\|\rot\Big(\bE_\alpha-\bE_0-\alpha\bw_0(\frac{\bx-\bz}{\alpha})\Big)
\Big\|_{L^2(B_\alpha)} \leq C \alpha^{7/2}(|1-\mu_r^{-1}|+ \nu)
\|\nabla\times\bE_0\|_{W^{2,\infty}(B_\alpha)},\\
\nm &
&\Big\|\bE_\alpha-\bE_0-\nabla\phi_0-\alpha\bw_0\Big(\frac{\bx-\bz}{\alpha}\Big)\Big\|_{L^2(B_\alpha)}
\leq C \alpha^{9/2}(|1-\mu_r^{-1}|+ \nu)
\|\nabla\times\bE_0\|_{W^{2,\infty}(B_\alpha)}. \een
\end{thm}

To conclude this section we remark that
\be
\label{F}
\alpha^{-1}\bF(\bz+\alpha\bm\xi)&=&\i\omega\mu_0\left(\frac 12\,\bH_0(\bz)\times\bm\xi+\frac \alpha3\,
{\bf D}\bH_0(\bz)\bm\xi\times\bm\xi\right)\nn\\
\nm &=&\i\omega\mu_0\left(\frac 12\,\sum^3_{i=1}\bH_0(\bz)_i{\itbf
e}_i\times\bm\xi+\frac \alpha3\,\sum^3_{i,j=1}{\bf
D}\bH_0(\bz)_{ij}{\itbf e}_i{\itbf
e}_j^T\bm\xi\times\bm\xi\right), \ee where ${\bf
D}\bH_0(\bz)_{ij}$ is the $(i,j)$-th element of the matrix ${\bf
D}\bH_0(\bz)$ and $T$ denotes the transpose. Thus \be\label{w0}
\bw_0({\bm \xi})=\i\omega\mu_0 \left(\frac
12\sum^3_{i=1}\bH_0(\bz)_i{\bm \theta}_i({\bm \xi})+\frac
\alpha3\,\sum^3_{i,j=1}{\bf
D}\bH_0(\bz)_{ij}\bm\Psi_{ij}(\bm\xi)\right), \ee where ${\bm
\theta}_i({\bm \xi})$ is the solution of the following interface
problem
\begin{equation} \label{thetaH}
\left\{
\begin{array}{l}
\nabla_{\bm \xi} \times \mu^{-1} \nabla_{\bm \xi} \times{\bm \theta}_i-\i \omega\sigma\alpha^2 {\bm
\theta}_i= \displaystyle\i\omega\sigma\alpha^2{\itbf
e}_i\times\bm\xi\mbox{ in }B\cup B^c,\\\nm
 \nabla_{\bm \xi} \cdot {\bm
\theta}_i = 0 \mbox{ in } B^c,\\\nm \mbox{[}{\bm \theta}_i
\times\bn]_{\Gamma} =0,\ \ [\mu^{-1} \nabla_{\bm \xi} \times{\bm \theta}_i\times\bn
]_{\Gamma}= - 2 [\mu^{-1}]_\Gamma {\itbf e}_i\times \bn \mbox{ on } {\Gamma}, \\
\nm {\bm
\theta}_i({\bm \xi})=O(|\bm\xi|^{-1})\mbox{ as } |{\bm \xi}|
\rightarrow\infty,
\end{array}
\right.
\end{equation}
and ${\bm\Psi}_{ij}({\bm \xi})$ is the solution of
\begin{equation} \label{thetaE}
\left\{
\begin{array}{l}
\nabla_{\bm \xi} \times\mu^{-1} \nabla_{\bm \xi} \times{\bm\Psi}_{ij}-\i\omega\sigma\alpha^2{\bm\Psi}_{ij}=
\i\omega\sigma\alpha^2\xi_j{\itbf e}_i\times\bm\xi \mbox{ in
}B\cup B^c,\\\nm\nabla_{\bm \xi}  \cdot {\bm\Psi}_{ij} =0 \mbox{ in }
B^c,\\\nm \mbox{[}{\bm\Psi}_{ij} \times\bn]_{\Gamma} =0,\ \
[\mu^{-1} \nabla_{\bm \xi} \times{\bm\Psi}_{ij}\times\bn
]_{\Gamma}= - 3 [\mu^{-1}]_\Gamma  \xi_j {\itbf e}_i\times \bn \mbox{ on } {\Gamma}, \\
\nm {\bm\Psi}_{ij}({\bm \xi})=O(|\bm\xi|^{-1})\mbox{ as } |{\bm
\xi}| \rightarrow\infty.
\end{array}
\right.
\end{equation}
Here ${\itbf e}_i$ is the unit vector in the $x_i$ direction. It
is worth emphasizing that since $\nu = O(1)$, ${\bm \theta}_i$ and
${\bm\Psi}_{ij}$ are uniformly bounded in $\bX_1(\RR^3)$.

We impose $\nabla\cdot\bm\theta_i=0$ outside $B$ to make the
solution $\bm\theta_i$ unique outside $B$. In this case by
\cite[Proposition 3.1]{ABN} we can show that
$\bm\theta_i=O(|\bm\xi|^{-2})$ and
$\nabla\times\bm\theta_i=O(|\bm\xi|^{-3})$ as $|\bm\xi|\to\infty$,
which implies by integrating \eqref{thetaH} over $B$ that \ben
\i\omega\sigma_*\alpha^2\int_B(\bm\theta_{i}+ {\itbf{e}_i}\times\bm\xi)d\bm\xi&=&
\int_{\partial B}\bn\times\mu^{-1}\nabla\times\bm\theta_{i}d\bm\xi\\
\nm &=&\int_{\partial B_R}\bn\times\mu^{-1}\nabla\times\bm\theta_{i}d\bm\xi\\
\nm &\to&0\ \ \mbox{as } R \to +\infty, \een where $B_R$ is a ball
of radius $R$ so that $B\subset B_R$. Thus we obtain
\begin{equation}\label{thetaEB2} \int_B(\bm\theta_i+{\itbf
e}_i\times\bm\xi)d\,\bm\xi=0.
\end{equation}

Similarly, by imposing $\nabla\cdot\bm\Psi_{ij}=0$ outside $B$ we
know that $\nabla\times\bm\Psi_{ij}=O(|\bm\xi|^{-3})$. Moreover,
integrating \eqref{thetaE} over $B$ and using similar argument
leading to \eqref{thetaEB2} together with the symmetry of ${\bf
D}\bH_0(\bz)$ yields
 \be \label{thetaEB}
\ds \sum_{i,j=1}^3 {\bf D}\bH_0(\bz)_{ij}
\int_B(\bm\Psi_{ij}+\xi_j{\itbf{e}_i}\times\bm\xi)d \bm\xi=0. \ee

\subsection{Integral Representation Formulas}

The integral representation is similar to the Stratton-Chu formula
for time-harmonic Maxwell equations (cf., e.g., N\'ed\'elec
\cite{nedelec}).

\begin{lem}\label{L:3.3} Let $D$ be a bounded domain in $\RR^3$ with Lipschitz boundary $\Gamma_D$ whose unit outer normal is $\bn$. For any $\bE\in\mathbf{H}_{-1}(\mathbf{curl};\RR^3\backslash\bar D)$
satisfying $\rot\rot\bE=0, \ddiv\bE=0$ in $\RR^3\backslash\bar D$, we have, for any $\bx\in \RR^3\backslash\bar D$,
\begin{eqnarray*}
\bE(\bx)&=&-\nabla_\bx \times \int_{\Gamma_D}(\bE(\by)\times\bn)G(\bx,\by)d \by
-\int_{\Gamma_D}(\nabla_\by \times \bE(\by)\times\bn)G(\bx,\by) d \by \\
&&
-\nabla_\bx\int_{\Gamma_D}(\bE(\by)\cdot\bn)G(\bx,\by)d \by,
\end{eqnarray*}
 where $G(\bx,\by)=\frac 1{4\pi |\bx-\by|}$ is the fundamental
solution of the Laplace equation.
\end{lem}

\debproof For the sake of completeness we give a sketch of proof.
Since $\bE\in\mathbf{H}_{-1}(\mathbf{curl};\RR^3\backslash\bar
D)$, for any $\bF$ such that $\bF(\by)=O(|\by|^{-1})$ and ${\bf D}
\bF (\by)=O(|\by|^{-2})$ as $|\by|\to\infty$, we can obtain by
integrating by parts, the conditions
$\nabla\times\nabla\times\bE=0, \nabla\cdot\bE=0$ in
$\RR^3\backslash\bar D$,  that \ben
(\bE,-\Delta\bF)_{\RR^3\backslash\bar D}
&=&(\bE,\rot\rot\bF-\nabla\ddiv\bF)_{\RR^3\backslash\bar D}\ \ \ \ \ \\
\nm
&=&-\int_{\Gamma_D}(\bE\times\bn)\cdot\rot\bF d \by
-\int_{\Gamma_D}\rot\bE\times\bn\cdot\bF d\by
+\int_{\Gamma_D}(\bE\cdot\bn)\ddiv\bF d\by.
\end{eqnarray*}
Now for $\bx\in\RR^3\backslash\bar D$ and $j\in \{1,2,3\}$, we
choose $\bF(\by)=G(\bx,\by){\itbf e}_j$ and thus $-\Delta_\by
\bF=\delta(\bx,\by){\itbf e}_j$, where $\delta(\bx,\cdot)$ is the
Dirac mass at $\bx$. Then we have
\begin{eqnarray*}
\nm
\bE_j(\bx)
&=&-\int_{\Gamma_D}(\bE(\by)\times\bn)\cdot\nabla_\by\times(G(\bx,\by){\itbf e}_j)d \by
-\int_{\Gamma_D}(\nabla_\by \times \bE(\by)\times\bn)_jG(\bx,\by) d\by \\
&&+\int_{\Gamma_D}(\bE(\by)\cdot\bn)\frac{\partial G(\bx,\by)}{\partial y_j} d\by\\
\nm
&=&-\left(\nabla_\bx \times \int_{\Gamma_D}(\bE(\by)  \times\bn)G(\bx,\by) d\by \right)_j
-
\int_{\Gamma_D}(\nabla_\by\times \bE (\by) \times\bn)_jG(\bx,\by) d\by
\\
&&
-\frac{\partial}{\partial x_j}\int_{\Gamma_D}(\bE(\by) \cdot\bn)G(\bx,\by) d\by,
\end{eqnarray*}
where we have used the fact that
\begin{eqnarray*}
(\bE(\by)\times\bn)\cdot\nabla_\bx\times(G(\bx,\by){\itbf e}_j)=-(\nabla_\bx\times(G(\bx,\by)\bE(\by)\times\bn))_j.
\end{eqnarray*}
This completes the proof. \finproof

The following lemma will be useful in deriving the asymptotic formula in next subsection.
Recall that $\displaystyle \bH_\alpha=\frac{1}{\i\omega\mu_\alpha}\rot\bE_\alpha, \bH_0=\frac{1}{\i\omega\mu_0}\rot\bE_0$.

\begin{lem}\label{L:3.4}
Let $\tilde\bH_\alpha=\bH_\alpha-\bH_0$. Then we have, for $\bx\in B_\alpha^c$,
 \begin{eqnarray*}
 \tilde\bH_\alpha(\bx)=\int_{B_\alpha}\nabla_\bx G(\bx,\by)
 \times\nabla_\by \times \tilde\bH_\alpha(\by)\,d\by+(1-\frac{\mu_*}{\mu_0})
 \int_{B_\alpha}(\bH_\alpha(\by)\cdot\nabla_\by)\nabla_\bx G(\bx,\by)\,d\by.
 \end{eqnarray*}
\end{lem}

\debproof
It is easy to check that $\rot\tilde\bH_\alpha=0$ and $\ddiv\tilde\bH_\alpha=0$ in $B^c_\alpha$.
By the representation formula in Lemma \ref{L:3.3} we have
$$
\tilde\bH_\alpha(\bx)=-\nabla_\bx \times \int_{\Gamma_\alpha}(\tilde\bH^+_\alpha(\by)\times\bn)G(\bx,\by)d\by
-\nabla_\bx\int_{\Gamma_\alpha}(\tilde\bH^+_\alpha(\by)\cdot\bn)G(\bx,\by) d\by,
$$
where $\tilde\bH_\alpha^+=\tilde\bH_\alpha|_{B_\alpha^c}$. Denote
$\tilde\bH_\alpha^-=\tilde\bH_\alpha|_{B_\alpha}$ and let
$\bE^{\pm}_\alpha$ be defined likewise. By the interface condition
$[\tilde\bE_\alpha\times\bn]_{\Gamma_\alpha}=0$, we have
\begin{eqnarray*}
\tilde\bH_\alpha^+\cdot\bn=\frac{1}{\i\omega\mu_0}\nabla
\times\bE_\alpha^+\cdot\bn-\bH_0\cdot\bn&=&
\frac{1}{\i\omega\mu_0}\rm{div}_{\Gamma_\alpha}(\bE^+_\alpha\times\bn)-\bH_0\cdot\bn\\
\nm &=&\frac{\mu_*}{\mu_0}\bH^-_\alpha\cdot\bn-\bH_0\cdot\bn,
\end{eqnarray*}
where $\rm{div}_{\Gamma_\alpha}$ denotes the surface divergence.
Then since $[\tilde\bH_\alpha\times\bn]_{\Gamma_\alpha}=0$, we
have
\begin{eqnarray}
\nonumber
\tilde\bH_\alpha(\bx)
&=&-\nabla_\bx \times \int_{\Gamma_\alpha}(\tilde\bH^-_\alpha(\by)\times\bn)G(\bx,\by)d\by\\
&&-\nabla_\bx \int_{\Gamma_\alpha}(\frac{\mu_*}{\mu_0}\bH^-_\alpha(\by)\cdot\bn-\bH_0(\by)\cdot\bn)G(\bx,\by)d\by.
\label{h-alpha}
\end{eqnarray}
For the first term,
\begin{eqnarray}\label{x5}
&&-\nabla_\bx \times \int_{\Gamma_\alpha}(\tilde\bH_\alpha^-(\by)\times\bn)G(\bx,\by) d\by\nn\\
\nm
&&=\nabla_\bx\int_{B_\alpha}\nabla_\by\times(\tilde\bH_\alpha(\by)G(\bx,\by)) d\by\nn\\
\nm
&&=\nabla_\bx\int_{B_\alpha}(G(\bx,\by)\nabla_\by \times \tilde\bH_\alpha(\by)+\nabla_\by G(\bx,\by)\times\tilde\bH_\alpha(\by)) d\by\nn\\
\nm
&&=
\int_{B_\alpha}\left(\nabla_\bx
G(\bx,\by)\times\nabla_\by \times \tilde\bH_\alpha(\by)
+(\tilde\bH_\alpha\cdot\nabla_\bx)\nabla_\by
G(\bx,\by)\right) d\by,
\end{eqnarray}
where we have used the identity
\ben
\rot(\bu\times\bv)=\bu(\ddiv\bv)-(\bu\cdot\nabla)\bv+(\bv\cdot\nabla)\bu-\bv(\ddiv\bu),
\een
and the fact that $\nabla_\bx\cdot\nabla_\by G(\bx,\by)=-\Delta_\by G(\bx,\by)=0$.
For the second term, we first notice that
\begin{eqnarray*}
-\nabla_\bx \int_{\Gamma_\alpha}(\frac{\mu_*}{\mu_0}\bH^-_\alpha(\by)\cdot\bn-\bH_0(\by)\cdot\bn)G(\bx,\by)d\by
&=&-\frac{\mu_*}{\mu_0}\nabla_\bx\int_{\Gamma_\alpha}\tilde\bH^-_\alpha(\by)\cdot\bn G(\bx,\by)d \by\\
&& +(1-\frac{\mu_*}{\mu_0}) \nablaÐ\bx \int_{\Gamma_\alpha}\bH_0(\by)\cdot\bn G(\bx,\by) d \by.
\end{eqnarray*}
By integration by parts we have
\begin{eqnarray*}
 \nabla_\bx\int_{\Gamma_\alpha}\tilde\bH_\alpha^-(\by)\cdot\bn G(\bx,\by) d\by
 &=&\nabla_\bx\int_{B_\alpha}\nabla_\by\cdot(G(\bx,\by)\tilde\bH_\alpha(\by)) d\by\\
 \nm
 &=&\nabla_\bx\int_{B_\alpha}\nabla_\by G(\bx,\by)\cdot\tilde\bH_\alpha(\by) +G(\bx,\by)\ddiv\tilde\bH_\alpha(\by) d\by\\
 \nm
 &=&\int_{B_\alpha}(\tilde\bH_\alpha(\by)\cdot\nabla_\by)\nabla_\bx G(\bx,\by) d\by.
 \end{eqnarray*}
 Similarly
 $$
 \nabla_\bx
 \int_{\Gamma_\alpha}(\bH_0(\by)\cdot\bn)G(\bx,\by)d\by
 =\int_{B_\alpha}(\bH_0(\by)\cdot\nabla_\by)\nabla_\bx G(\bx,\by) d\by.
 $$
Thus
 \begin{eqnarray}\label{x6}
& &-\nabla_\bx\int_{\Gamma_\alpha}(\frac{\mu_*}{\mu_0}\bH^-_\alpha(\by)\cdot\bn-\bH_0(\by)\cdot\bn)G(\bx,\by) d \by\nn\\
\nm
 &&=-\frac{\mu_*}{\mu_0}\int_{B_\alpha}(\tilde\bH_\alpha(\by)\cdot\nabla_\by)\nabla_\bx G(\bx,\by) d\by\nn\\
 \nm
 && \quad +\ (1-\frac{\mu_*}{\mu_0})\int_{B_\alpha}(\bH_0(\by)\cdot\nabla_\by)\nabla_\bx G(\bx,\by) d\by.
 \end{eqnarray}
This completes the proof by substituting \eqref{x5}-\eqref{x6} into \eqref{h-alpha}.
\finproof

\subsection{Asymptotic Formulas}

In this subsection we prove the following theorem which is the main result of this section.

\begin{thm}\label{thm-2} Let $\nu$ be of order one and let $\alpha$ be small. For $\bx$ away from the location $\bz$ of the inclusion, we have
\ben 
\bH_\alpha(\bx)-\bH_0(\bx)&=&\i\nu\alpha^3\left[\frac 12\,
\sum^3_{i=1}\bH_0(\bz)_i\int_B{\bm D}^2_\bx G(\bx,\bz)\bm\xi\times(\bm\theta_i
+{\itbf e}_i\times\bm\xi)d\bm\xi\right]\nn\\
\nm &&+\alpha^3\Big(1-\frac{\mu_0}{\mu_*}\Big)
\left[\sum^3_{i=1}\bH_0(\bz)_i {\bm D}^2_\bx G(\bx,\bz)
\int_B\Big({\itbf e}_i+\frac 12\nabla\times\bm\theta_i\Big)
d\bm\xi\right]+\bR(\bx), \een where $({\bm D}^2_\bx G )_{ij}=
\partial_{x_i x_j}^2 G$ and  \ben |\bR(\bx)|\le
C \alpha^4\|\bH_0\|_{W^{2,\infty}(B_\alpha)}, \een uniformly in
$\bx$ in any compact set away from $\bz$.

\end{thm}

\debproof The proof starts from the integral representation
formula in Lemma \ref{L:3.4}. We first consider the first term in
the integral representation in Lemma \ref{L:3.4}. By Theorem
\ref{T:3.1} we know that \be\label{y3}
\|\bE_\alpha-\bE_0-\nabla\phi_0-\alpha\bw_0(\frac{\bx-\bz}{\alpha})\|_{L^2(B_\alpha)}
\leq C\alpha^{9/2}(|1-\mu_r^{-1}|+ \nu)\|\nabla \times
\bE_0\|_{W^{2,\infty}(B_\alpha)}. \ee Since $\rot \bH_0=0 $ and
$\rot\bH_\alpha=\sigma\bE_\alpha$ in $B_\alpha$, we have
$$
\int_{B_\alpha}\nabla_\bx G(\bx,\by)\times\nabla_\by \times \tilde\bH_\alpha(\by) d\by = \sigma_*\int_{B_\alpha}\nabla_\bx G(\bx,\by)\times\bE_\alpha(\by) d\by
=
{\rm I_1}+\cdots+{\rm I_4},
$$
where
\begin{eqnarray*}
\nm
{\rm I_1}
&=&
\sigma_*\int_{B_\alpha}\nabla_\bx G(\bx,\by)\times\Big(\bE_\alpha(\by)
-\bE_0(\by)-\nabla\phi_0(\by)-\alpha\bw_0(\frac{\by-\bz}{\alpha})\Big) d\by , \\
{\rm I_2}
&=&
\sigma_*\int_{B_\alpha}\nabla_\bx G(\bx,\by)\times(\bE_0(\by) +\nabla\phi_0(\by)-\bF(\by))d\by , \\
{\rm I_3}
&=&
\sigma_*\int_{B_\alpha}(\nabla_\bx G(\bx,\by)-\nabla_\bx G(\bx,\bz)-{\itbf D}^2_{\bx}G(\bx,\bz)(\by-\bz))\times\Big(\bF(\by)+\alpha\bw_0(\frac{\by-\bz}\alpha)\Big)d\by , \\
{\rm I_4}
&=&
\sigma_*\int_{B_\alpha}(\nabla_\bx G(\bx,\bz)+{\itbf D}^2_{\bx}G(\bx,\bz)(\by-\bz))\times\Big(\bF(\by)+\alpha\bw_0(\frac{\by-\bz}{\alpha})\Big) d\by  .
\end{eqnarray*}
By \eqref{y3}, we have \ben |\,{\rm I}_1|&\le&C
\alpha^6(|1-\mu_r^{-1}|+ \nu)\sigma_*
\|\nabla\times\bE_0\|_{W^{2,\infty}(B_\alpha)}\\
\nm
&\le&Ck\alpha^6|1-\mu_r^{-1}|\|\bH_0\|_{W^{2,\infty}(B_\alpha)}\\
\nm &\le&C \alpha^4\|\bH_0\|_{W^{2,\infty}(B_\alpha)}.
 \een By \eqref{F} we have $|\,{\rm I}_2|\le
C\alpha^{6}\sigma_*\|\nabla\times\bE_0\|_{W^{2,\infty}(B_\alpha)}\le
C\alpha^4\|\bH_0\|_{W^{2,\infty}(B_\alpha)}$. Similarly, by using
\eqref{deff} and \eqref{w0} we can show $|{\rm I}_3|\le
C\alpha^4\|\bH_0\|_{W^{2,\infty}(B_\alpha)}$. For the remaining
term we first observe that \ben {\rm I}_4=\i
\alpha^4\sigma_*\int_B(\nabla_\bx G(\bx,\bz)+\alpha {\itbf
D}^2_{\bx}G(\bx,\bz)\bm\xi)\times(\alpha^{-1}\bF(z+\alpha\bm\xi)+\bw_0(\bm\xi))
d\bm\xi. \een On the other hand, \ben
\alpha^{-1}\bF(\bz+\alpha\bm\xi)+\bw_0(\bm\xi)&=&\i\omega\mu_0\Big[\frac
12\,\sum^3_{i=1}\bH_0(\bz)_i({\itbf
e}_i\times\bm\xi+\bm\theta_i({\bm \xi})
)\\
&&+\frac \alpha3\,\sum^3_{i,j=1}{\bm D}\bH_0(\bz)_{ij}(\xi_j{\bm e}_i\times\bm\xi+\bm \Psi_{ij}({ \bm \xi}))\Big],
\een
which implies after
using  \eqref{thetaEB}
\ben
{\rm I}_4&=&\i k\alpha^4\left[\frac 12\,\sum^3_{i=1}\bH_0(\bz)_i\int_B\nabla_\bx G(\bx,\bz)\times({\itbf e}_i\times\bm\xi+\bm\theta_i)d\,\bm\xi\right]\\
\nm &&+\i k\alpha^5\left[\frac
12\,\sum^3_{i=1}\bH_0(\bz)_i\int_B{\bm D}^2_\bx
G(\bx,\bz)\bm\xi\times({\itbf
e}_i\times\bm\xi+\bm\theta_i)d\,\bm\xi\right]
 + \bR_1(\bx),
\een where $|\bR_1(\bx)|\le
C\alpha^4\|\bH_0\|_{W^{2,\infty}(B_\alpha)}$. Using
\eqref{thetaEB2}, this shows that \be \label{y4}
& &\int_{B_\alpha}\nabla_\bx G(\bx,\by)\times\nabla_\by \times\tilde\bH_{\alpha}(\by)d\,\by\nn\\
\nm &=& \i k\alpha^5\left[\frac
12\,\sum^3_{i=1}\bH_0(\bz)_i\int_B{\bm D}^2_\bx
G(\bx,\bz)\bm\xi\times({\itbf
e}_i\times\bm\xi+\bm\theta_i)d\,\bm\xi\right] +\bR_2(\bx), \ee
where $$|\bR_2(\bx)|\le C
k\alpha^6|1-\mu_r^{-1}|\|\bH_0\|_{W^{2,\infty}(B_\alpha)} \le C
\alpha^4\|\bH_0\|_{W^{2,\infty}(B_\alpha)}.$$

Now we turn to the second term in Lemma \ref{L:3.4}. From Theorem
\ref{T:3.1} we know that \be\label{y1}
\left\|\bH_\alpha-\frac{\mu_0}{\mu_*}\bH_0-\frac{\alpha}{\i\omega\mu_*}\nabla_\bx\times
\bw_0(\frac{\bx-\bz}{\alpha})\right\|_{L^2(B_\alpha)}\leq
C\alpha^{7/2} (|1-\mu_r^{-1}|+ \nu)
\|\bH_0\|_{W^{2,\infty}(B_\alpha)}. \ee
 Let
$$
\displaystyle\bH_0^*({\bm
\xi})=\frac{1}{\i\omega\mu_0}\nabla_{\bm \xi}\times\bw_0({\bm
\xi}) .
$$
 Then
$$
 \int_{B_\alpha}(\bH_\alpha\cdot\nabla_\by)\nabla_\bx G(\bx,\by) d\by
 = -\int_{B_\alpha}{\bm D}^2_\bx G(\bx,\by)\bH_\alpha(\by)d\by
= {\rm II}_1+\cdots+{\rm II}_4  ,
 $$
 where
 \begin{eqnarray*}
{\rm II}_1  &=&-\int_{B_\alpha}{\bm D}^2_\bx G(\bx,\by)\left(\bH_\alpha(\by)-\frac{\mu_0}{\mu_*}\bH_0(\by)
 -\frac{\mu_0}{\mu_*}\bH^*_0(\frac{\by-\bz}{\alpha})\right)d\by,\\
{\rm II}_2&=&-\frac{\mu_0}{\mu_*}\int_{B_\alpha}({\bm D}^2_\bx G(\bx,\by)-{\bm D}^2_\bx(\bx,\bz))(\bH_0(\by)+\bH^*_0(\frac{\by-\bz}{\alpha}))d\by , \\
{\rm II}_3 &=&-\frac{\mu_0}{\mu_*}\int_{B_\alpha}{\bm D}^2_\bx G(\bx,\bz)(\bH_0(\by)-\bH_0(\bz))d\by, \\
{\rm II}_4 &=&-\frac{\mu_0}{\mu_*}\int_{B_\alpha}{\bm D}^2_\bx G(\bx,\bz)(\bH_0(\bz)+\bH^*_0(\frac{\by-\bz}{\alpha})) d\by.
 \end{eqnarray*}
 It is easy to see from \eqref{y1} that $|\,{\rm II}_1|\le C\alpha^4\|\bH_0\|_{W^{1,\infty}(B_\alpha)}$.
 By \eqref{w0} we know that
 \ben
 \| \bH^*_0(\frac{\by-\bz}{\alpha})\|_{L^2(B_\alpha)}\le C\alpha^{3/2} \|\bH_0\|_{W^{2,\infty}(B_\alpha)},
 \een
which implies $|\,{\rm II}_2|\le C\alpha^4\|\bH_0\|_{W^{2,\infty}(B_\alpha)}$.
Similarly, we have
$|\,{\rm II}_3|\le C\alpha^4 \|\bH_0\|_{W^{1,\infty}(B_\alpha)}$.
Finally, by \eqref{w0}, we have
\ben {\rm II}_4=-\frac{\mu_0}{\mu_*}\alpha^3\sum^3_{i=1}\bH_0(\bz)_i\int_B {\bm D}^2_\bx
G(\bx,\bz)\Big({\itbf e}_i+\frac 12\nabla\times\bm\theta_i\Big)d\,\bm\xi+\bR_3(\bx),
\een
where $|\bR_3(\bx)|\le C\alpha^4\|\bH_0\|_{W^{2,\infty}(B_\alpha)}$. Therefore,
\be
\label{y2}
& &\int_{B_\alpha}(\bH_\alpha\cdot\nabla_\by)\nabla_\bx G(\bx,\by)
\,d\by\nn\\
\nm &=& -\frac{\mu_0}{\mu_*}\alpha^3\sum^3_{i=1}\bH_0(\bz)_i\int_B
{\bm D}^2_\bx G(\bx,\bz)\Big({\itbf e}_i+\frac
12\nabla\times\bm\theta_i\Big)d\,\bm\xi+ \bR_4(\bx)
\ee
with $|\bR_4(\bx)|\le
C\alpha^4 \|\bH_0\|_{W^{2,\infty}(B_\alpha)}$. This completes the
proof by substituting \eqref{y2} and \eqref{y4} into the integral
representation formula in Lemma \ref{L:3.4}. \finproof

It is worth emphasizing that the tensor whose column vectors are
$$\Big(1-\frac{\mu_0}{\mu_*}\Big) \int_B \Big({\itbf e}_1+\frac
12\nabla\times\bm\theta_1\Big) d\bm\xi,
\Big(1-\frac{\mu_0}{\mu_*}\Big) \int_B \Big({\itbf e}_2+\frac
12\nabla\times\bm\theta_2\Big) d\bm\xi,\quad \mbox{and }
\Big(1-\frac{\mu_0}{\mu_*}\Big) \int_B \Big({\itbf e}_3+\frac
12\nabla\times\bm\theta_3\Big) d\bm\xi$$ is the so-called magnetic
polarization tensor. It reduces in the zero conductivity case
($\sigma=0$) to the one first introduced in \cite{AVV}.

On the other hand, for an arbitrary shaped target, one introduces
for $l,l^\prime=1,2,3$, $\mathbb{M}^{(l,l^\prime)}$ to  be the
$3\times 3$ matrix whose $i$-th column is $$\ds \frac 12
 {\itbf e}_l   \times \int_B \xi_{l^\prime} (\bm\theta_i+{\itbf
e}_i\times\bm\xi)d\bm\xi.$$ One can easily show that
\begin{equation} \label{arbitrayshaped}
\ds \frac{1}{2} \sum^3_{i=1} \bH_0(\bz)_i \int_B{\bm D}^2_\bx
G(\bx,\bz)\bm\xi\times(\bm\theta_i+{\itbf e}_i\times\bm\xi)d\bm\xi
 =\ds  \sum^3_{l,l^\prime=1} {\bm D}^2_\bx G(\bx,\bz)_{ll^\prime} \mathbb{M}^{(l,l^\prime)} \bH_0(\bz). \end{equation}
 We call $\mathbb{M}^{(l,l^\prime)}$ the conductivity
polarization tensors.

We now consider the case of a spherical target. If $B$ is a
sphere, then one can check that
$$
\frac 12\sum^3_{i=1}\bH_0(\bz)_i \int_B{\bm D}^2_\bx
G(\bx,\bz)\bm\xi\times(\bm\theta_i+{\itbf e}_i\times\bm\xi)d\bm\xi
=  \mathcal{M} {\bm D}^2_\bx G(\bx,\bz) \bH_0(\bz),
$$
where \begin{equation} \label{defm} \mathcal{M} = \frac{1}{2}
\int_B (\xi_1 \bm\theta_2(\bm\xi) \cdot {\itbf e}_3 - \xi_1^2)\,
d\bm\xi,\end{equation} and therefore, the asymptotic formula
derived in Theorem \ref{thm-2} reduces in the case $\mu_0=\mu_*$
to the following result.

\begin{cor} \label{corM}
Assume that $\mu_0=\mu_*$ and $B$ is a sphere. Then we have \be
\label{sigmaonly2} \bH_\alpha(\bx)-\bH_0(\bx) &=& \i k\alpha^5
\mathcal{M} {\bm D}_{\bx}^2G(\bm\bx,\bm\bz) \bH_0(\bz) + \bR(\bx).
\ee  The remainder satisfies $|\bR(\bx)|\le
C\alpha^4\|\bH_0\|_{W^{2,\infty}(B_\alpha)}$ uniformly in $\bx$ in
any compact set away from $\bz$.
\end{cor}

Now we assume that $\bJ_0$ is a dipole source whose position is
denoted by $\bs$
\begin{equation}
 \bJ_0(\bx)= \rot \big( \bp\,
\delta(\bx, \bs) \big)  ,
\end{equation}
where $\delta(\cdot, \bs)$ is the Dirac  mass at $\bs$ and the
unit vector $\bp$ is the direction of the magnetic dipole. The
existence and uniqueness of a solution to \eqref{alpha} follows
from \cite{valli}. In the absence of any inclusion, the magnetic
field $\bH_0$ due to $\bJ_0(\bx)$ is given by
\begin{equation}
\label{ehon} \bH_0(\bx) = \rot \rot (\bp G(\bx, \bs))={\bm D}^2_{\bx}G(\bx,\bs)\bp.
\end{equation}

We note that $J_0$ is not in the dual of $\bX_\alpha(\RR^3)$,
however we can form the difference $\bE_\alpha - \bE_0$ and solve
for that difference in $\widetilde{\bX}_\alpha(\RR^3)$. That way
we are able to recover Theorems \ref{T:3.1} and \ref{thm-2}.

The asymptotic formula \eqref{sigmaonly2} can be rewritten as
\begin{equation} \label{finalformula}
\bq \cdot (\bH_\alpha-\bH_0)(\bx) \simeq  \i k \alpha^5
\mathcal{M} \left({\bm D}^2_\bx G(\bx,\bz)\bq\right)^T \left({\bm
D}^2_\bx G(\bz,\bs)\bp\right),
\end{equation}
where $m$ is defined by \eqref{defm}. Note that, in view of
\eqref{finalformula}, if the dipole $\bJ_0$ is located at $\bx$,
then the field $\bp\cdot \bH_\alpha$ at $\bs$ is the same as the
one obtained if $\bJ_0$ is located at $\bs$ and $\bp\cdot
\bH_\alpha$ is measured at $\bx$.

Next, writing
$$
\mathcal{M}  = \Re e \, \mathcal{M}  + \i \Im m \,  \mathcal{M} ,
$$
we obtain
$$
\Re e \big( \bq \cdot (\bH_\alpha-\bH_0)(\bx) \big) \simeq - k
\alpha^5 (\Im m \,\mathcal{M} ) \left({\bm D}^2_\bx
G(\bx,\bz)\bq\right)^T \left({\bm D}^2_\bx G(\bz,\bs)\bp\right),
$$
and
$$
\Im m \big( \bq \cdot (\bH_\alpha-\bH_0)(\bx) \big) \simeq  k
\alpha^5 (\Re e \,\mathcal{M} ) \left({\bm D}^2_\bx
G(\bx,\bz)\bq\right)^T \left({\bm D}^2_\bx G(\bz,\bs)\bp\right).
$$
In view of \eqref{ehon}, $\bH_0$ is real. Therefore, it follows
that
\begin{equation}
\label{realfinal} \Im m \big(  \bq \cdot \bH_\alpha(\bx) \big)
\simeq  k \alpha^5 (\Re e \,\mathcal{M} ) \left({\bm D}^2_\bx
G(\bx,\bz)\bq\right)^T \left({\bm D}^2_\bx G(\bz,\bs)\bp\right).
\end{equation}
Formula (\ref{realfinal}) will be used in the section for locating
and detecting a spherical target. For arbitrary shaped targets,
the formula derived in Theorem \ref{thm-2} together with
(\ref{arbitrayshaped}) should be used.

\section{Localization and Characterization}
In this section we consider that there are ${M}$ sources and $N$ receivers.
The $m$th source  is located at $\bs_m$ and it
 generates the magnetic dipole ${\itbf J}_0^{(m)} ({\itbf r})
=\nabla \times ({\itbf p} \delta( {\itbf r}, \bs_m) )$.
The $n$th receiver is located at $\br_n$ and it records the magnetic
field in the ${\itbf q}$ direction.
The $(n,m)$-th entry of the $N\times M$ response matrix ${\bf A}=(A_{nm})_{n=1,\ldots,N,m=1,\ldots, {M}}$
is the signal recorded by the $n$th receiver when the $m$th source is emitting:
$$
A_{nm}=
 \bH_\alpha^{(m)}(\br_n) \cdot \bq .
$$
The response matrix is the sum of the unperturbed field $
\bH_0^{(m)}(\br_n) \cdot \bq $ and the perturbation $
\bH_\alpha^{(m)}(\br_n) \cdot \bq - \bH_0^{(m)}(\br_n) \cdot \bq
$. This perturbation contains information about the inclusion but
it is much smaller (of order $\alpha^3$) than the unperturbed
field. Consequently, it seems that we need to know the unperturbed
field with great accuracy in order to be able to extract the
perturbation and to process it. In practice, such an accuracy may
not be accessible. However, we know that the unperturbed field is
real while the perturbation is complex-valued, as shown by
(\ref{realfinal}). The imaginary part of the response matrix is
therefore equal to the imaginary part of the perturbation and this
is the data that we will process:
\begin{equation}
\label{eq:unpert}
{\bf A}_0=(A_{0,nm})_{n=1,\ldots,N,m=1,\ldots, {M}}, \quad \quad
{A}_{0,nm}= \Im m \,\big(\bH_\alpha^{(m)}(\br_n) \cdot
\bq\big) .
\end{equation}

We assume that $N \geq M$, that is, there are more receivers than
sources. As in \cite{AILP}, in order to locate the conductive
inclusion $\bz+\alpha B$ we can use the MUSIC imaging functional.
We focus on formula (\ref{realfinal})
 and  define the MUSIC imaging functional for a search point $\bz^S$ by
\begin{equation} \label{musicfunct}
{\mathcal I}_{\rm MU}(\bz^S)= \bigg[\frac{1}{\sum_{l=1}^3 \|
\left({\bf I}_N - {\bf P}  \right)({\bm D}^2_\bx G(\br_1,\bz^S)\bq
\cdot \bev_l, \ldots, {\bm D}^2_\bx G(\br_N,\bz^S)\bq \cdot
\bev_l)^T \|^2}\bigg]^{1/2},
\end{equation}
where ${\bf P}$ is the orthogonal projection on the range of the
matrix ${\bf A}_0$ and $(\bev_1,\bev_2,\bev_3)$ is an orthonormal
basis of $\RR^3$. From \cite{AILP}, it follows that the following
proposition holds.
\begin{prop}
In the presence of an inclusion located at $z$, the matrix ${\bf
A}_0$ has three significant singular values counted with their
multiplicity. Moreover, the MUSIC imaging functional ${\mathcal
I}_{\rm MU}$ attains its maximum approximately at $z^S=z$.
\end{prop}

Once the inclusion is located we can compute by a least-squares
method   $\Re e \, \mathcal{M} $ associated with the inclusion
from the response matrix ${\bf A}_0$. Given the location of the
inclusion, we minimize the discrepancy between the computed and
the measured response matrices. For a single frequency, knowing
$\Re e\, \mathcal{M} $ may not be sufficient to separate the
conductivity of an inclusion from its size. However, $\Re e\,
\mathcal{M} $ obtained for a few frequencies $\omega$ may be used
to reconstruct both the conductivity and the size of the target.

\section{Noisy Measurements}
In this section we consider that there are ${M}$ sources and $N$
receivers. The measures are noisy, which means that the magnetic
field measured by a receiver is corrupted by an additive noise
that can be described in terms of a real Gaussian random variable
with mean zero and variance $\sigma_{\rm n}^2$. The recorded
noises are independent from each other.

\subsection{Hadamard Technique}
{\bf Standard acquisition.} In the standard acquisition scheme,
the response matrix is measured at each step of ${M}$ consecutive
experiments. In the $m$th experience, $m=1,\ldots,{M}$, the $m$th
source (located at $\bs_m$)
 generates the magnetic dipole ${\itbf J}_0^{(m)} ({\itbf r})
=\nabla \times ({\itbf p} \delta( {\itbf r}, \bs_m) )$ and the $N$
receivers (located at $\br_n$, $n=1,\ldots,N$) record the magnetic
field in the ${\itbf q}$ direction
which means that they measure
$$
A_{{\rm meas},nm} = A_{0,nm} + W_{nm}, \quad n=1,\ldots,N, \quad m=1,\ldots,{M}  ,
$$
which gives the matrix
\begin{equation}
\label{Ameas1}
{\bf A}_{{\rm meas}} = {\bf A}_{0} + {\bf W} ,
\end{equation}
where ${\bf A}_0$ is the unperturbed response matrix (\ref{eq:unpert})
and $W_{nm}$ are independent Gaussian random variables with mean
zero and variance $\sigma_{\rm n}^2$. Here,
$\bH_\alpha^{(m)}(\br_n)$ is the magnetic field generated by a
magnetic dipole at $\bs_m$ and measured at the receiver $\br_n$ in
the presence of the inclusion.

 The so called Hadamard noise reduction technique is valid in the
 presence of additive noise and uses the structure of Hadamard matrices.

\begin{definition}
A Hadamard matrix ${\bf H}$ of order ${M}$  is a ${M} \times {M}$
matrix whose elements are $-1$ or $+1$ and such that ${\bf H}^T
{\bf H}= {M} {\bf I}_M$. Here ${\bf I}_M$ is the $M\times M$ identity
matrix.
\end{definition}
Hadamard matrices do not exist for all ${M}$.
A necessary condition for the existence is that ${M}=1,2$ or a multiple of $4$.
A sufficient condition is that  ${M}$ is a power of two.
Explicit examples are known for all ${M}$ multiple of $4$ up to ${M}=664$ \cite{seberry}.\\

\noindent {\bf Hadamard acquisition.} In the Hadamard acquisition
scheme, the response matrix is measured during a sequence of ${M}$
experiments. In the $m$th experience, $m=1,\ldots,{M}$, all
sources generate magnetic dipoles, the $m'$ source generating $
H_{mm'} {\itbf J}_0^{(m')} ({\itbf r})$. This means that we use
all sources to their maximal emission capacity (which is a physical constraint) with a specific
coding of their signs. The $N$ receivers record the magnetic field
in the ${\itbf q}$ direction,
which means that they measure
$$
B_{{\rm meas},nm} =  \sum_{m'=1}^{M} H_{mm'} A_{0,nm'} + W_{nm} = ({\bf A}_0 {\bf H}^T)_{nm}
+ W_{nm} , \quad n=1,\ldots,N,\quad m=1,\ldots,{M}  ,
$$
which gives the matrix
$$
{\bf B}_{\rm meas} ={\bf A}_0 {\bf H}^T +{\bf  W},
$$
where ${\bf A}_0$ is the unperturbed response matrix
and $W_{nm}$ are independent Gaussian random variables with mean zero and variance $\sigma_{\rm n}^2$.
The measured response matrix ${\bf A}_{\rm meas}$
is obtained by right multiplying
the matrix ${\bf B}_{\rm meas}$ by  the matrix $\frac{1}{{M}} {\bf H}$:
$$
{\bf A}_{{\rm meas}} = \frac{1}{{M}}  {\bf B}_{\rm meas} {\bf H}
=  \frac{1}{{M}}  {\bf A}_0 {\bf H}^T {\bf H} +  \frac{1}{{M}}
{\bf W}{\bf H}  ,
$$
which gives
\begin{equation}
{\bf A}_{{\rm meas}} = {\bf A}_0   + \widetilde{\bf W}, \quad \quad  \widetilde{\bf W}= \frac{1}{{M}} {\bf  W} {\bf H}.
\end{equation}
The benefit of using Hadamard's technique lies in the fact that
the new noise matrix $ \widetilde{\bf W}$ has independent entries
with Gaussian statistics, mean zero, and variance $\sigma_{\rm
n}^2 /{M}$:
\begin{eqnarray*}
\EE\big[  \widetilde{W}_{nm}  \widetilde{W}_{n'm'} \big] &=&
\frac{1}{{M}^2} \sum_{q,q'=1}^{M} H_{qm} H_{q'm'} \EE[ W_{nq}
W_{n'q'} ] =
\frac{\sigma_{\rm n}^2}{{M}^2} \sum_{q,q'=1}^{M} H_{qm} H_{q'm'} \delta_{nn'} \delta_{qq'}\\
&=& \frac{\sigma_{\rm n}^2}{{M}^2} \sum_{q=1}^{M} H_{qm}
(H^T)_{m'q} \delta_{nn'} =
\frac{\sigma_{\rm n}^2}{{M}^2}  ({\bf H}^T {\bf H})_{m'm} \delta_{nn'} \\
&=& \frac{\sigma_{\rm n}^2}{{M}}   \delta_{mm'}  \delta_{nn'} ,
\end{eqnarray*}
where $\EE$ stands for the expectation and $\delta_{mn}$ is the
Kronecker symbol. This gain of a factor ${M}$ in the
signal-to-noise ratio is called the Hadamard advantage.

\subsection{Singular Values of a Noisy Matrix}
We consider in this subsection the case where there is measurement
noise but no inclusion is present. The measured response matrix is
the $N\times M$ matrix
\begin{equation}
\label{eq:withoutinc}
{\bf A}_{\rm meas} = {\bf W},
\end{equation}
where ${\bf W}$ consists of independent noise coefficients
with mean zero and variance $\sigma_{\rm n}^2/{M}$ and the number of receivers
is larger than the number of sources $N \geq {M}$.
This is the case when the response matrix is acquired with the Hadamard technique and there is no
inclusion in the medium.

We denote by  $  \sigma_{1}^{({M})}  \geq \sigma_{2}^{({M})} \geq
\sigma_{3}^{({M})}  \geq \cdots \geq \sigma_{{M}}^{({M})}$ the singular
values of the  response matrix ${\bf A}_{\rm meas}$ sorted by decreasing order
and by $\Lambda^{({M})}$  the corresponding integrated density of states defined by
$$
\Lambda^{({M})}([a,b]) = \frac{1}{{M}}
{\rm Card} \left\{ l=1,\ldots,{M}  \,  , \, \sigma_l^{({M})}  \in [  a, b] \right\}, \quad
\mbox{for any } a<b  .
$$
The density $\Lambda^{({M})}$ is a counting measure which consists
of a sum of Dirac masses:
$$
\Lambda^{({M})} = \frac{1}{{M}} \sum_{j=1}^{M} \delta_{ \sigma_j^{({M})}  }   .
$$
For large $N$ and ${M}$ denote $$ \gamma:= N / {M} \geq 1.$$ The
following statements hold.
\begin{prop}
\label{prop:dens:2}%
\begin{enumerate}
\item[{\rm a)}] The random measure $\Lambda^{({M})}$  converges
almost surely to the deterministic absolutely continuous measure
$\Lambda$ with compact support:
\begin{equation}
\label{semicircle1}
\Lambda ([\sigma_u,\sigma_v]) = \int_{\sigma_u}^{\sigma_v}
\frac{1}{\sigma_{\rm n}}
\rho_{\gamma} \Big(\frac{\sigma}{\sigma_{\rm n}}\Big)
d\sigma,
\quad \quad 0\leq \sigma_u \leq \sigma_v ,
\end{equation}
where
$\rho_\gamma$ is the deformed quarter-circle law given by
\begin{equation}
\label{semicircle}
 \rho_\gamma (\sigma) = \left\{
\begin{array}{ll}
\displaystyle  \frac{1}{\pi \sigma }
 \sqrt{ \big( (\gamma^{1/2}+1)^2 - \sigma^2 \big) \big( \sigma^2- (\gamma^{1/2}-1)^2\big)} &
\mbox{ if } \gamma^{1/2}-1 < \sigma \leq \gamma^{1/2}+1 , \\
\displaystyle 0 & \mbox{ otherwise} .
\end{array}
\right.
\end{equation}
\item[{\rm b)}]
The normalized $l^2$-norm of the singular values satisfies
\begin{equation}
{M} \Big[
\frac{1}{{M}} \sum_{j=1}^{M}  ( \sigma^{({M})}_j )^2 - \gamma \sigma_{\rm n}^2  \Big]
 \stackrel{{M} \to \infty}{\longrightarrow}  \sqrt{2 \gamma} \sigma_{\rm n}^2 Z
  \mbox{ in distribution},
 \end{equation}
where $Z$ follows a Gaussian distribution with mean zero and variance one.
\item[{\rm c)}]
The maximal singular value satisfies
\begin{equation}
\label{maxsigma1}
\sigma_1^{({M})} \simeq \sigma_{\rm n} \Big[
\gamma^{1/2} + 1  + \frac{1}{2{M}^{2/3}}  \big( 1+\gamma^{-1/2}\big)^{1/3}  Z_1 + o( \frac{1}{{M}^{2/3}}) \Big]
  \mbox{ in distribution},
\end{equation}
where $Z_1$ follows a type-1 Tracy-Widom distribution.
\end{enumerate}
\end{prop}


The type-1 Tracy-Widom distribution has the pdf $p_{\rm TW1}$:
\begin{eqnarray*}
&&\PP(Z_1 \leq z)  = \int_{-\infty}^z p_{{\rm TW}1}(x) dx= \exp
\Big(
 -\frac{1}{2} \int_z^\infty \varphi(x) +(x-z) \varphi^2(x) dx \Big)  ,
 \end{eqnarray*}
where $\varphi$ is the solution of the Painlev\'e
equation
\begin{equation} \label{painv} \varphi''(x)=x \varphi(x) +2\varphi(x)^3, \quad \quad \varphi(x)
\stackrel{x \to +\infty}{\simeq} {\rm Ai}(x),
\end{equation} ${\rm Ai}$
being the Airy function.
The expectation of $Z_1$ is $\EE[ Z_1]  \simeq -1.21$ and its variance is $ {\rm Var}(Z_1) \simeq
1.61$.

\debproof
 Point a) is Marcenko-Pastur result \cite{pastur2}.  Point b)
follows from the expression of the normalized $l^2$-norm of the
singular values in terms of the entries of the matrix:
$$
\frac{1}{{M}} \sum_{j=1}^{M}  ( \sigma^{({M})}_j )^2 =
\frac{1}{{M}} {\rm tr} \big( {\bf A}_{\rm meas}^T {\bf A}_{\rm meas} \big) =
\frac{1}{{M}} \sum_{n=1}^N \sum_{m=1}^{M} A_{{\rm meas},nm}^2,
$$
and from the application of the central limit theorem in the regime ${M} \gg 1$.
The third point follows from \cite{johnstone01}.
\finproof

\subsection{Singular Values of the Unperturbed Response Matrix}
We now turn to the case where there is one conductive inclusion in
the medium and no measurement noise. The measured response matrix
  is then
 the $N \times {M}$ matrix ${\bf A}_0$ defined by
\begin{equation}
 \label{A0}
A_{0,nm} =   k \alpha^5 (\Re \, e \mathcal{M} ) \big({\bm D}^2_\bx
G(\br_n,\bz)\bq \big)^T \,  \,  \big({\bm D}^2_\bx G(\bz,\bs_m)\bq
\big) .
\end{equation}
The matrix ${\bf A}_0 $ possesses three nonzero singular values
given by
\begin{eqnarray*}
\sigma^{{\bf A}_0}_j &=& k \alpha^5  |\Re \, e \mathcal{M}|
\Big[\sum_{m=1}^{M}
 \Big| ( \big( {\bm D}^2_\bx G(\bz,\bs_m)\bq\big) \big)_j
\Big|^2 \Big]^{1/2} \\
&& \times \Big[\sum_{n=1}^N
 \Big| ( \big(  {\bm D}^2_\bx G(\br_n,\bz)\bq  \big) \big)_j
\Big|^2 \Big]^{1/2} , \quad j=1,2,3 .
\end{eqnarray*}


\subsection{Singular Values of the Perturbed Response Matrix}
The measured
 response matrix using the Hadamard technique
 in the presence of an inclusion and in the presence of measurement noise is
 \begin{equation}
 \label{AW}
{\bf A}_{\rm meas} = {\bf A}_0 +{\bf W}  ,
\end{equation}
where ${\bf A}_0$ is given by (\ref{A0}) and ${\bf W}$
has independent random entries
with Gaussian statistics,
mean zero and variance $\sigma_{\rm n}^2/{M}$.

We consider the critical regime in which the singular values of
the unperturbed matrix are of the same order as the singular
values of the noise, that is to say, $\sigma_1^{{\bf A}_0}$, the
first singular value of ${\bf A}_0 $,  is of the same order
of magnitude as $\sigma_{\rm n}$. We will say a few words about
the cases $\sigma_1^{{\bf A}_0}$ much larger or much smaller than
$\sigma_{\rm n}$ after the analysis of the critical regime.

\begin{prop}
\label{prop:dens:3}%
\begin{enumerate}
\item[{\rm a)}]
The normalized $l^2$-norm of the singular values satisfies
\begin{equation}
{M} \Big[
\frac{1}{{M}} \sum_{j=1}^{M}  ( \sigma^{({M})}_j )^2 - \gamma \sigma_{\rm n}^2  \Big]
 \stackrel{{M} \to \infty}{\longrightarrow} (\sigma_0^{{\bf A}_0})^2 + \sqrt{2 \gamma} \sigma_{\rm n}^2 Z
  \mbox{ in distribution},
 \end{equation}
where $Z$ follows a Gaussian distribution with mean zero and variance one and
\begin{equation}
\sigma_0^{{\bf A}_0}  = \Big[ \sum_{j=1}^3 (\sigma_j^{{\bf A}_0} )^2 \Big]^{1/2}  .
\end{equation}
\item[{\rm b1)}] If $\sigma_1^{{\bf A}_0} <    \gamma^{1/4}
\sigma_{\rm n}$, then the maximal singular value satisfies
\begin{equation}
\sigma_1^{({M})} \simeq \sigma_{\rm n} \Big[
\gamma^{1/2} + 1  + \frac{1}{2{M}^{2/3}}  \big( 1+\gamma^{-1/2}\big)^{1/3}  Z_1 + o( \frac{1}{{M}^{2/3}}) \Big]
  \mbox{ in distribution},
\end{equation}
where $Z_1$ follows a type-1 Tracy Widom distribution. \item[{\rm
b2)}] If $\sigma_1^{{\bf A}_0} >  \gamma^{1/4} \sigma_{\rm n}$,
then
\begin{equation}
\sigma_1^{({M})} = \sigma_1^{{\bf A}_0} \Big( \alpha + O
(\frac{1}{{M}^{1/2}})  \Big) \mbox{ in probability},
\end{equation}
where
\begin{equation}
\label{def:alpha}
\alpha = \bigg(1+ (1+ \gamma) \frac{\sigma_{\rm n}^2}{
(\sigma_1^{{\bf A}_0} )^2} + \gamma \frac{\sigma_{\rm n}^4}{
(\sigma_1^{{\bf A}_0} )^4 }\bigg)^{1/2}.
\end{equation}
 If, additionally, $\sigma_1^{{\bf A}_0} > \sigma_2^{{\bf
A}_0} $, then the maximal singular value in the regime ${M} \gg 1$
has Gaussian distribution with the mean and variance given by
\begin{eqnarray}
\EE \big[ \sigma_1^{({M})} \big]
&=&
\sigma_1^{{\bf A}_0}
 \Big( \alpha
 + o (\frac{1}{{M}^{1/2}})\Big)
,\\
{\rm Var}\big( \sigma_1^{({M})} \big) &=&\frac{\sigma_{\rm
n}^2}{{M}} \Big( \beta + o (1) \Big),
\end{eqnarray}
where
\begin{equation}
\label{def:beta}
\beta= \frac{1- \gamma \frac{\sigma_{\rm n}^4}{ (\sigma_1^{{\bf
A}_0} )^4 }}{\bigg( 1+ (1+\gamma) \frac{\sigma_{\rm n}^2}{
(\sigma_1^{{\bf A}_0} )^2 } + \gamma \frac{\sigma_{\rm n}^4}{
(\sigma_1^{{\bf A}_0} )^4 }\bigg)^{1/2}} .
\end{equation}
\end{enumerate}
\end{prop}

\debproof Point a) follows again from the explicit expression of
the $l^2$-norm of the singular values in terms of the entries of
the matrix. Point b) in the case $N={M}$ is addressed in
\cite{capitaine09b} and the extension to $N \geq {M}$ is only
technical and can be obtained from
\cite{benaych,garniersolnamsri}. \finproof

Note that, in the item b2), if $\sigma_1^{{\bf A}_0} =
\sigma_2^{{\bf A}_0} \geq \sigma_3^{{\bf A}_0} $, then the
fluctuations are not Gaussian any more, but they can be
characterized as shown in~\cite{capitaine09b}. Note also that
formula (\ref{def:beta}) seems to  predict
 that the variance of the maximal singular value cancels
when $\sigma_1^{{\bf A}_0}  \searrow   \gamma^{1/4 }  \sigma_{\rm n} $, but this is true only to
the order $M^{-1}$, and in fact it becomes of order $M^{-4/3}$.
Following \cite{baik05} we can anticipate that there are interpolating distributions which appear when
$\sigma_1^{{\bf A}_0}=   \gamma^{ {1}/{4}}   \sigma_{\rm n} +  {w}{M^{-{1}/{3}}}  $
for some fixed $w$.

\subsection{Detection Test}
%

The objective in this subsection is to design  a detection method
which comes with an estimate of the level of
  confidence, in the presence of noise, in our ability to determine whether
  there actually is a conductive inclusion.

Since we know that the presence of an inclusion is characterized
by the existence of three significant singular values for ${\bf
A}_0 $, we propose to use a test of the form $R
> r$ for the alarm corresponding to the presence of a conductive
inclusion. Here $R$ is the quantity obtained from the measured
response matrix defined by
\begin{equation}
\label{def:Rk}
R = \frac{ \sigma^{({M})}_1 }{\Big[ \frac{1}{{M}-3(1+\gamma^{-1/2})^2}
\sum_{j=4}^{M} ( \sigma_j^{({M})} )^2 \Big]^{1/2}} ,
\end{equation}
and the threshold value $r$ has to be chosen by the user.
This choice follows from Neyman-Pearson theory as we explain below.
It requires the knowledge of the statistical distribution of $R$ which we give
in the following proposition.

\begin{prop}
In the asymptotic regime ${M} \gg 1$ the following statements hold.
\begin{enumerate}
\item[{\rm a)}]
In absence of a conductive inclusion (Eq.~(\ref{eq:withoutinc})) we have
\begin{equation}
\label{R3:dist1}
R \simeq
1 + \gamma^{-1/2}+ \frac{1}{2{M}^{2/3}}  \gamma^{-1/2} \big( 1+\gamma^{-1/2}\big)^{1/3}  Z_1 + o( \frac{1}{{M}^{2/3}})  ,
\end{equation}
where $Z_1$ follows a type-1 Tracy Widom distribution.
\item[{\rm b)}]
In presence of a conductive inclusion (Eq.~(\ref{AW})) :
\begin{enumerate}
\item[{\rm b1)}] If $\sigma_1^{{\bf A}_0} >  \gamma^{1/4}
\sigma_{\rm n}$, then we have
\begin{equation}
\label{R3:dist2} R \simeq \frac{\sigma_1^{{\bf A}_0}}{\gamma^{1/2}
\sigma_{\rm n}} \alpha + \frac{\beta^{1/2}}{\gamma^{1/2} {M}^{1/2}}
 Z_0  ,
\end{equation}
where $Z_0$ follows a Gaussian distribution with mean zero and
variance one. \item[{\rm b2)}] If $\sigma_1^{{\bf A}_0} <
\gamma^{1/4} \sigma_{\rm n}$, then we have (\ref{R3:dist1}).
\end{enumerate}
\end{enumerate}
\end{prop}

\debproof
In absence of a conductive inclusion, we have on the one hand
that the truncated normalized $l^2$-norm of the singular values satisfies
$$
{M} \Big[
\frac{1}{{M}-3(1+\gamma^{-1/2})^2} \sum_{j=4}^{M}  ( \sigma^{({M})}_j )^2 - \gamma \sigma_{\rm n}^2  \Big]
 \stackrel{{M} \to \infty}{\longrightarrow}
  \sqrt{2 \gamma} \sigma_{\rm n}^2 Z
  \mbox{ in distribution},
$$
where $Z$ follows a Gaussian distribution with mean zero and variance one,
which implies that
\begin{equation}
\label{truncl2}
\Big[ \frac{1}{{M}-3(1+\gamma^{-1/2})^2} \sum_{j=4}^{M} ( \sigma_j^{({M})} )^2 \Big]^{1/2}  =\gamma^{1/2} \sigma_{\rm n}
+ o (\frac{1}{{M}^{2/3}} )  \mbox{ in probability},
\end{equation}
and on the other hand the maximal singular value satisfies (\ref{maxsigma1}) in distribution.
Using Slutsky's theorem, we find the first item of the proposition.

In presence of a conductive inclusion, we have on the one hand
that the truncated normalized $l^2$-norm of the singular values satisfies (\ref{truncl2}).
On the other hand the maximal singular value is described by Proposition \ref{prop:dens:3}
which gives the desired result by Slutsky's theorem.
\finproof

%


%


The data (i.e. the measured response matrix) gives the value of the ratio
$R$. We propose to use a test of the form $R > r$ for the
alarm corresponding to the presence of a conductive inclusion.
The quality of this test can be quantified by two coefficients:\\
- The false alarm rate (FAR)  is the probability to sound the
alarm while there is no inclusion:
$$
{\rm FAR} =\PP( R > r | \mbox{ no inclusion} )  .
$$
-
The probability of detection (POD)  is   the probability to sound the alarm when there is an inclusion:
$$
{\rm POD}   = \PP ( R> r  | \mbox{ inclusion} )  .
$$

It is not possible to find a test that minimizes the FAR and maximizes the POD.
However, by the
Neyman-Pearson lemma, the decision rule of sounding the alarm
if and only if $R   > r_\delta$
maximizes the POD
 for a given FAR $\delta$ with the threshold
\begin{equation}
\label{ralpha}
r_\delta =   1+ \gamma^{-1/2} +    \frac{1}{2 {M}^{2/3}} \gamma^{-1/2}
\big(1 +\gamma^{-1/2}\big)^{1/3} \Phi_{{\rm TW}1} ^{-1} (1-\delta),
\end{equation}
 where $\Phi_{{\rm TW}1} $ is the cumulative distribution function of
the  Tracy-Widom distribution of type $1$. The computation of the
threshold $r_\delta$ is easy since it depends only on the number of sensors
$N$ and $M$ and on the  FAR  $\delta$.
Note that we should use a  Tracy-Widom
distribution table. We have, for
instance, $\Phi_{{\rm TW}1}^{-1}(0.9)\simeq 0.45$, $\Phi_{{\rm
TW}1}^{-1}(0.95)\simeq 0.98$ and $\Phi^{-1}_{{\rm TW}1}(0.99)
\simeq 2.02$.

The POD of this optimal test (optimal amongst all tests with the FAR $\delta$)
depends on the value $\sigma_1^{{\bf A}_0}$
and  on the noise level $\sigma_{\rm n}$.
Here we find that the POD  is
$$
{\rm POD}
 =\Phi \bigg(
\sqrt{{M}} \frac{ \frac{\sigma_1^{{\bf A}_0}}{  \sigma_{\rm n}}
\alpha -\gamma^{1/2} r_\delta }{\beta^{1/2}}\bigg),
$$
 where $\Phi $ is the cumulative distribution function of the normal distribution with mean zero and variance one.
The theoretical test performance   improves very rapidly with ${M}$
 once $\sigma_1^{{\bf A}_0} > \gamma^{1/4}
\sigma_{\rm n}$.
 This result is indeed valid as long as $\sigma_1^{{\bf A}_0} > \gamma^{1/4}
\sigma_{\rm n}$. When $\sigma_1^{{\bf A}_0} < \gamma^{1/4}
\sigma_{\rm n}$, so that the inclusion  is buried in noise (more
exactly, the singular values corresponding to the inclusion are
 buried into the deformed quarter-circle distribution of the other singular values),
 then we have
${\rm POD} =1 - \Phi_{{\rm TW}1} \big(  \Phi^{-1}_{{\rm TW}1}(1-\delta) \big) = \delta$.
 Therefore the probability of detection is given by
 \begin{equation}\label{pod}
{\rm POD} =   \max \bigg\{
  \Phi \bigg(
\sqrt{{M}} \frac{\frac{\sigma_1^{{\bf A}_0}}{  \sigma_{\rm n}}
\alpha -\gamma^{1/2} r_\delta }{\beta^{1/2}}\bigg) , \delta
\bigg\} .
 \end{equation}
 The transition region $\sigma_1^{{\bf A}_0} \simeq \gamma^{1/4}
\sigma_{\rm n}$ is only qualitatively characterized by our
analysis, as it would require a detailed study of the statistics
of the maximal singular value when $\sigma_1^{{\bf A}_0}=
\gamma^{ {1}/{4}}   \sigma_{\rm n} +  {w}{M^{-{1}/{3}}}  $ for
some fixed $w$.


Finally, the following remark is in order. The previous results
were obtained by an asymptotic analysis assuming that $M$ is large
and $\sigma_1^{{\bf A}_0}$ and $\sigma_{\rm n}$ are of the same
order. In the case in which $\sigma_1^{{\bf A}_0}$ is much larger
than $\sigma_{\rm n}$, then the proposed test has a POD of
$100\%$. In the case in which $\sigma_1^{{\bf A}_0}$ is much
smaller than $\sigma_{\rm n}$, then it is not possible to detect
the inclusion from the singular values of the response matrix and
the proposed test has a POD equal to the FAR (as shown above, this
is the case as soon as $\sigma_1^{{\bf A}_0} < \gamma^{1/4}
\sigma_{\rm n}$).

\section{Numerical Experiments}
In this section, we will give some numerical examples to
illustrate the performance of the detection algorithm. The
unperturbed measurement is acquired synthetically by asymptotic
formula \eqref{sigmaonly2} and noisy measurements are given by
\eqref{AW}. Assume that $B_\alpha$ is a sphere described by
$$(x-x_0)^2+(y-y_0)^2+(z-z_0)^2=\alpha^2,$$
where $\alpha$ is characteristic length of the inclusion measured
in meters. Then the domain $B$ is characterized by letting
$\alpha=1$ and $(x_0,y_0,z_0)$ to be origin. We assume that the
inclusion $B_\alpha$ is also located at the origin, $\alpha
=0.01$, $\mu_*=\mu_0=1.2566\times10^{-6}$~H/m and
$\sigma=5.96\times 10^7$~S/m. We let $\omega=133.5$ to make
$k\alpha^2=1$. We compute the solution of \eqref{thetaH} by an
edge element code. The numerically computed  $\mathcal{M}$ is
given by \begin{equation} \mathcal{M}=
   -0.4110 - 0.0387\i .
\end{equation}

The configuration of the detection system includes coincident
transmitter and receiver arrays uniformly distributed on the
square $[-2,2]\times[-2,2]\times \{ 1 \}$, both consisting of 256
($M=N=16^2$) vertical dipoles ($\bm p=\bm q=\bm e_3$) emitting or
receiving with unit amplitude. The search  domain is a box
$[-0.5, 0.5]^3$ below the arrays. It is worth mentioning here
that the number of transducers should be a multiple of $4$ in order
to be able to implement the Hadamard technique in a realistic situation.

\begin{figure}
\includegraphics[width=0.45\textwidth]{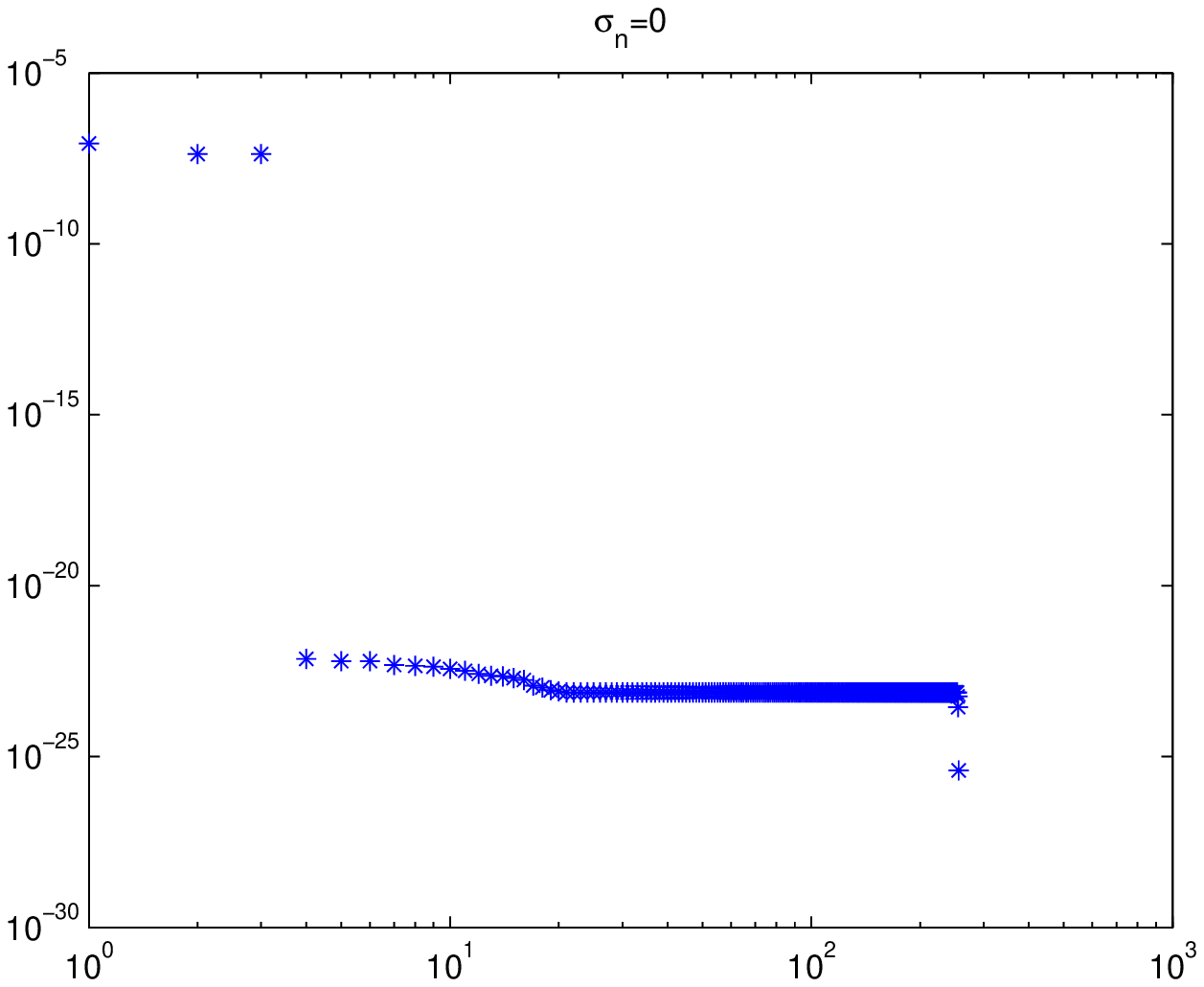}
\includegraphics[width=0.45\textwidth]{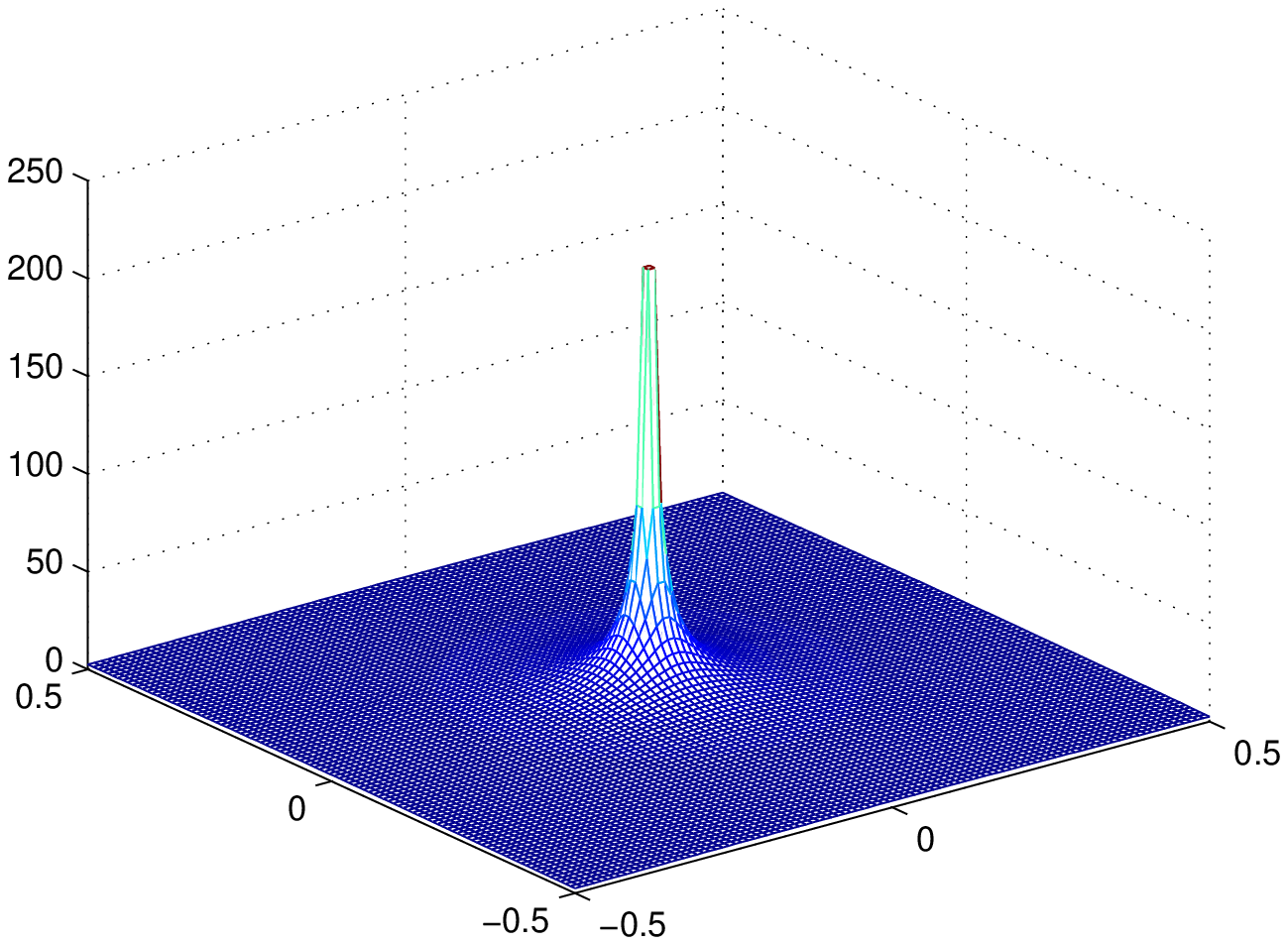}
\caption{Distribution of singular values of ${\bf A}_0$ with $M=N=256$ and the magnitude of $\mathcal{I}_{\rm MU}$ on plane $z=0.0$. } \label{fig6-1}
\end{figure}

 In the above setting, we calculate the SVD of the unperturbed response matrix ${\bf A}_0$. Figure \ref{fig6-1}
 displays the logarithmic scale plot of the singular values of ${\bf A}_0$.
We observe that our numerical results agree with our
 previous theoretical analysis: there is a significant singular value with multiplicity three associated with the inclusion.
 Then we can construct
 the projection ${\bf P}$  with the first three singular vectors corresponding to the first
 three significant singular values. In the right part of Figure \ref{fig6-1}, we also plot the magnitude
 of $\mathcal{I}_{\rm MU}$ on the cross section $z=0$, which shows that the MUSIC algorithm can detect the
 inclusion with high resolution.

\begin{figure}
\includegraphics[width=0.32\textwidth]{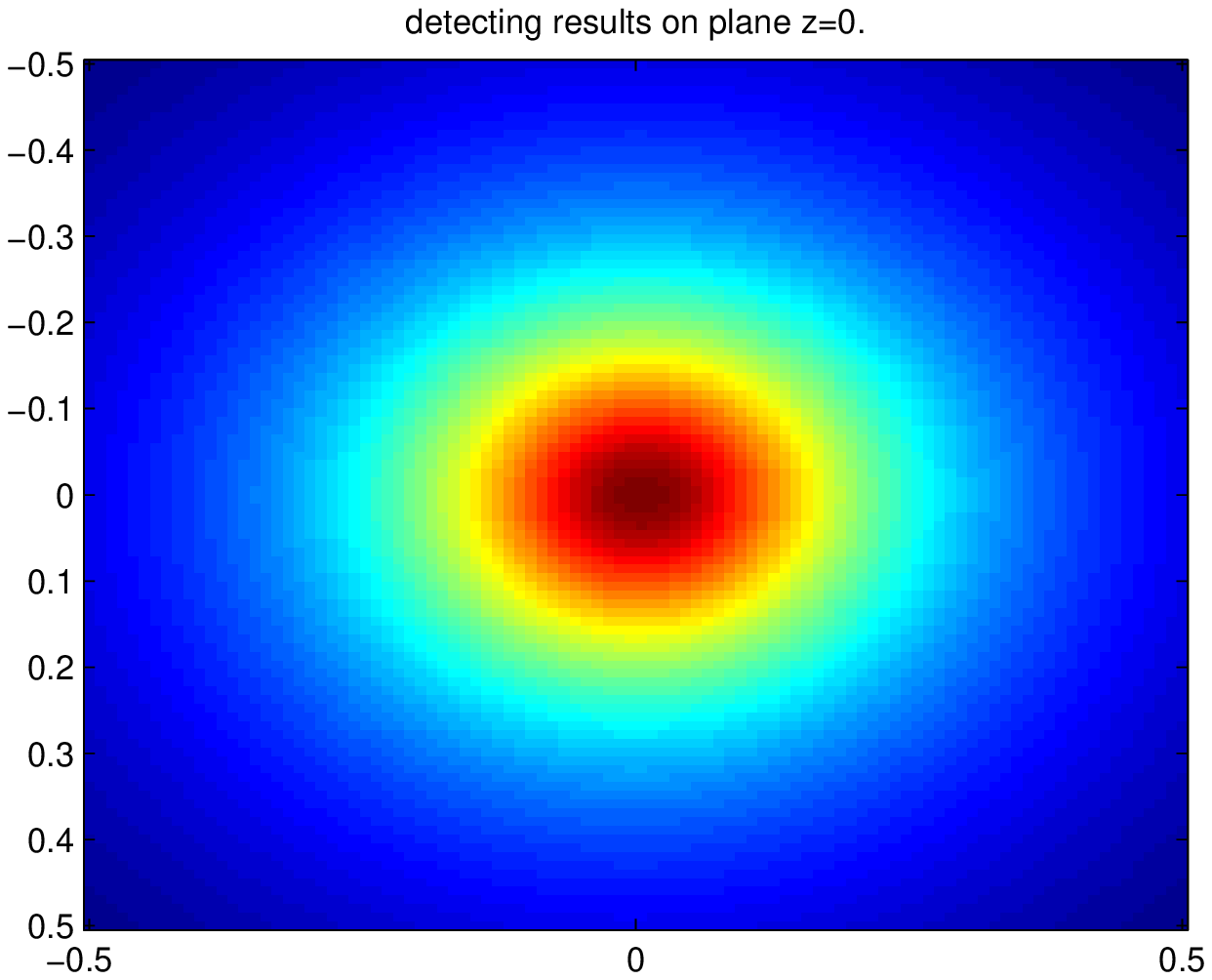}
\includegraphics[width=0.32\textwidth]{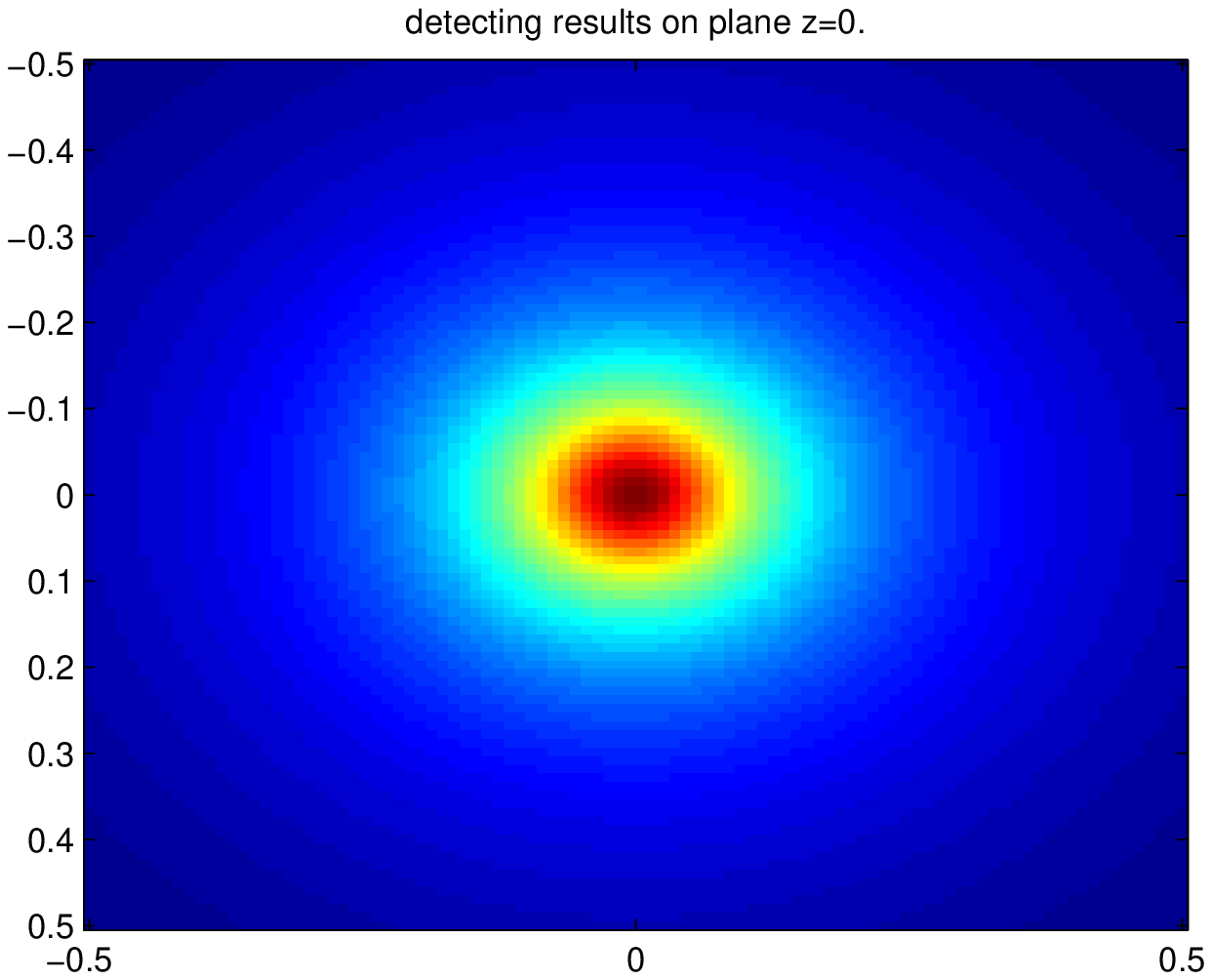}
\includegraphics[width=0.32\textwidth]{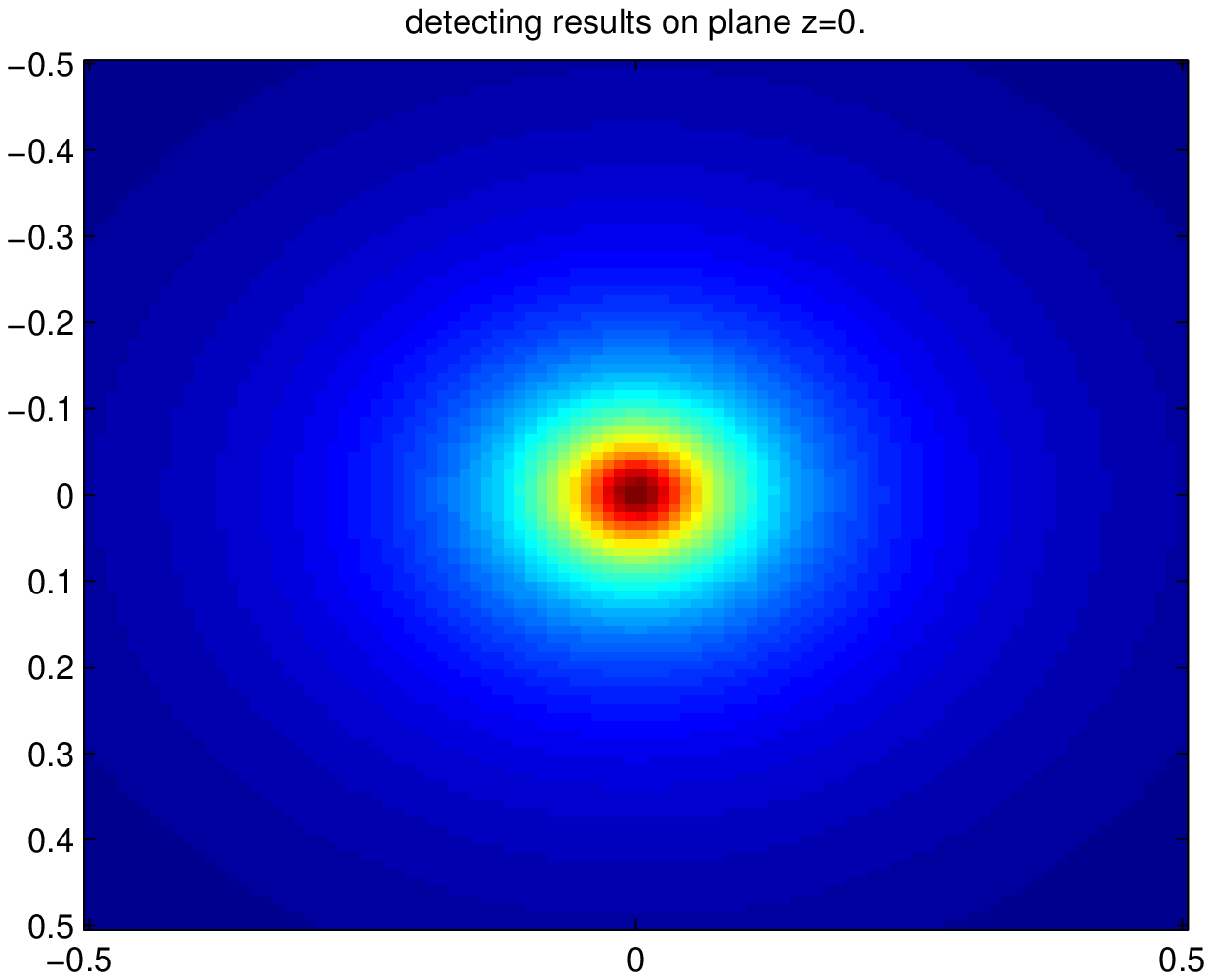}
\includegraphics[width=0.32\textwidth]{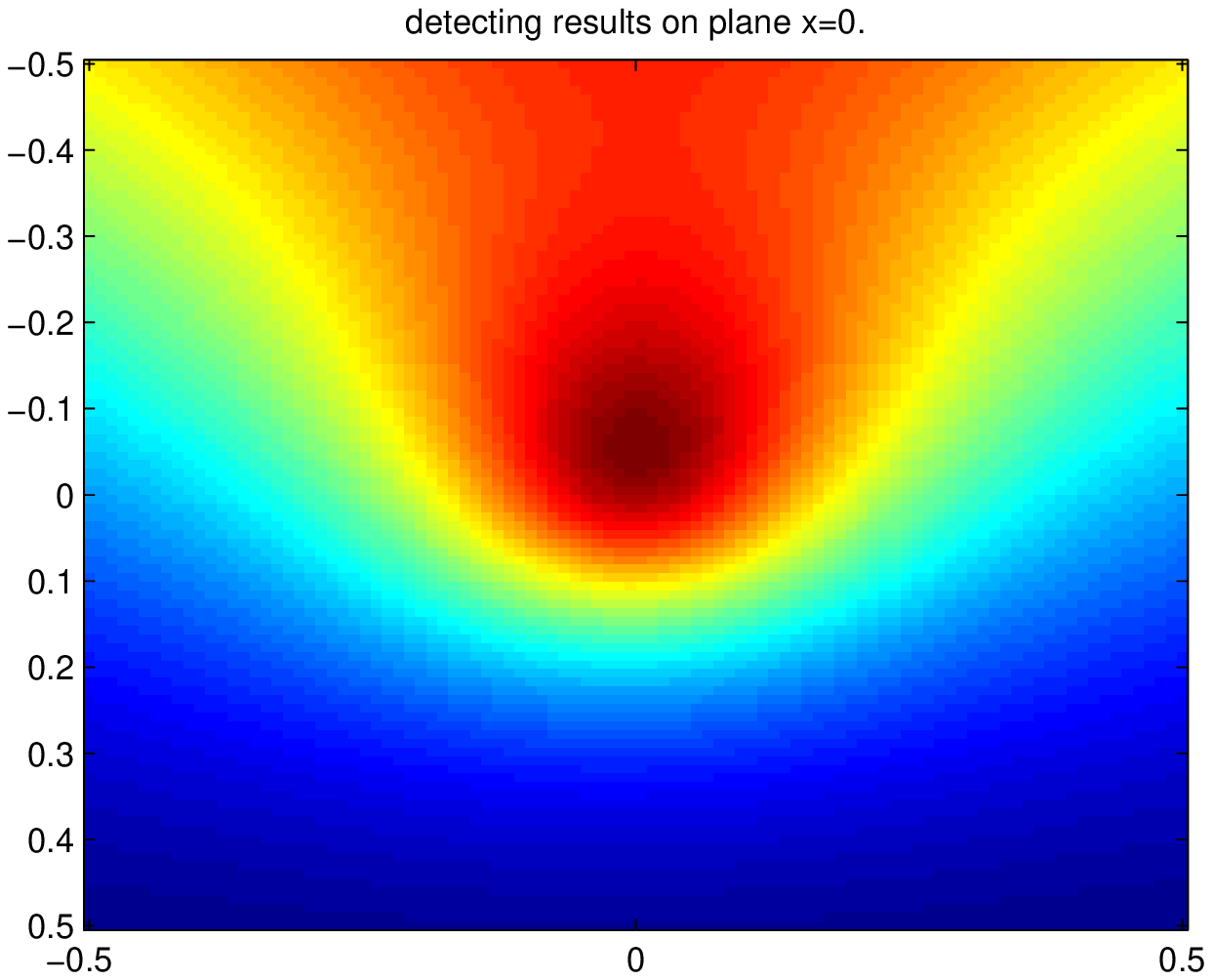}
\includegraphics[width=0.32\textwidth]{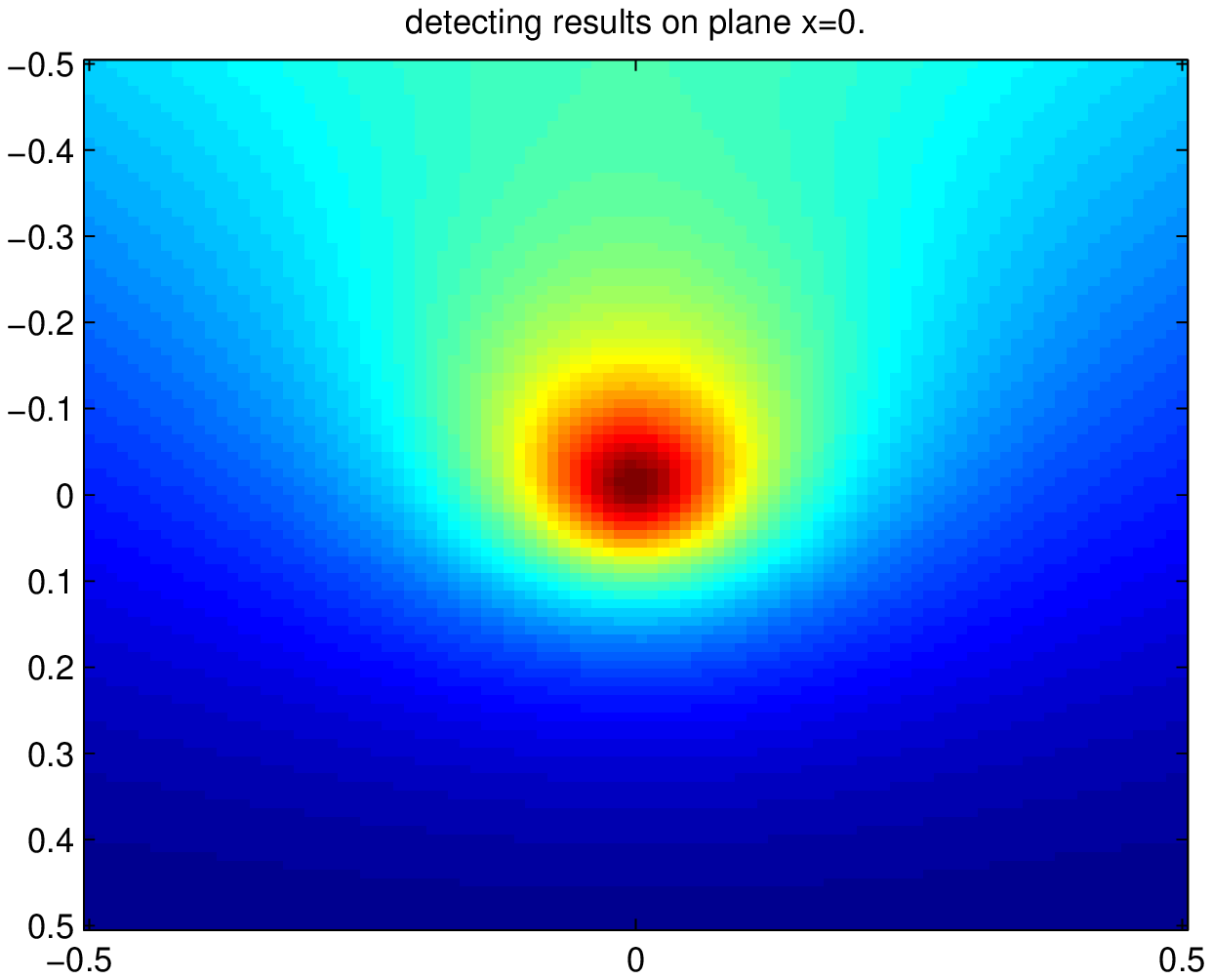}
\includegraphics[width=0.32\textwidth]{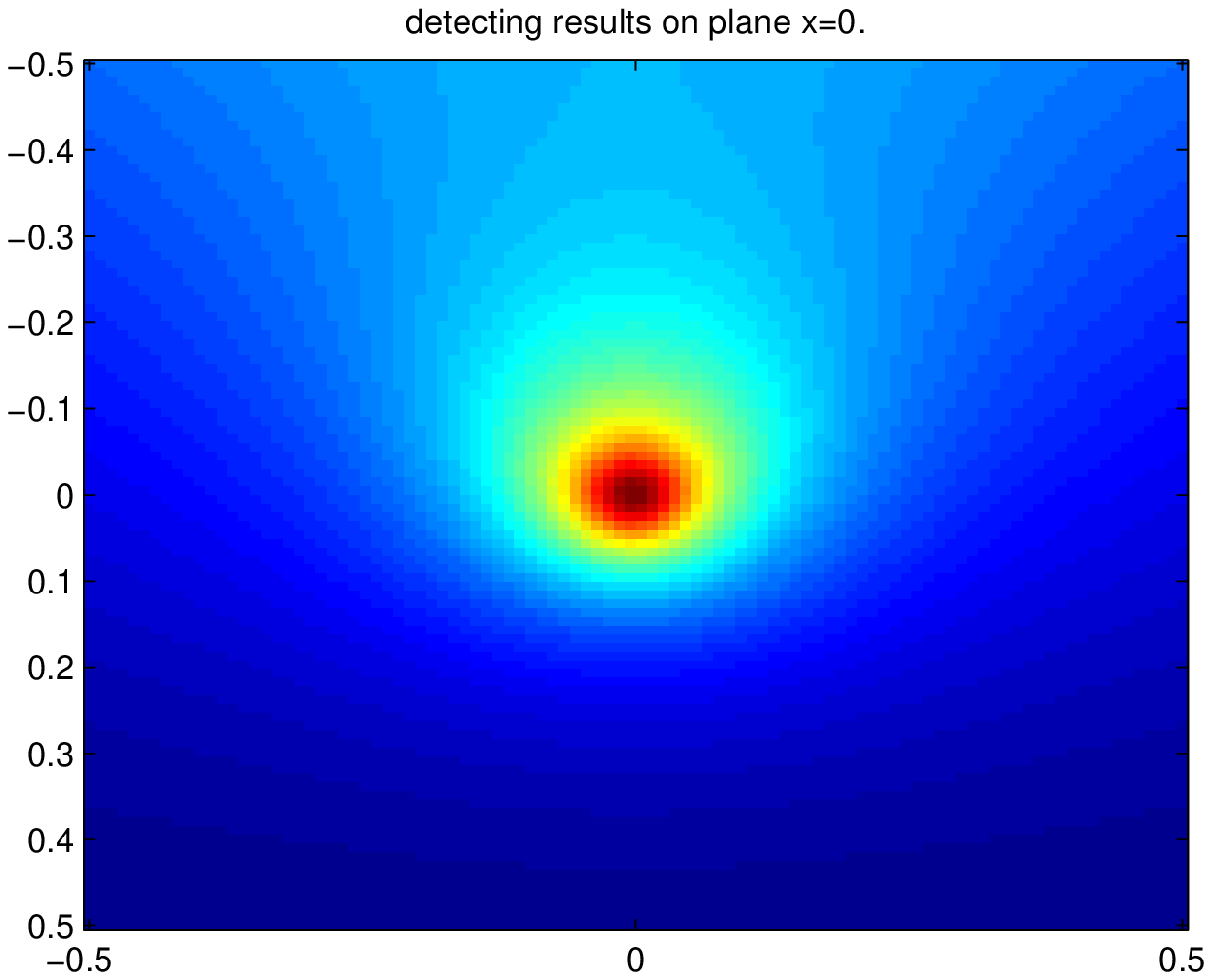}
\caption{Detecting results on cross sectional plane $z=0$(top)
and $x=0$(bottom) for different noise level $\sigma_{\rm n}$. $\sigma^{\bf A_0}_1/\sigma_{\rm n}=10,20,
30$ from left to right. }\label{fig6-2}
\end{figure}

\begin{figure}
\includegraphics[width=0.32\textwidth]{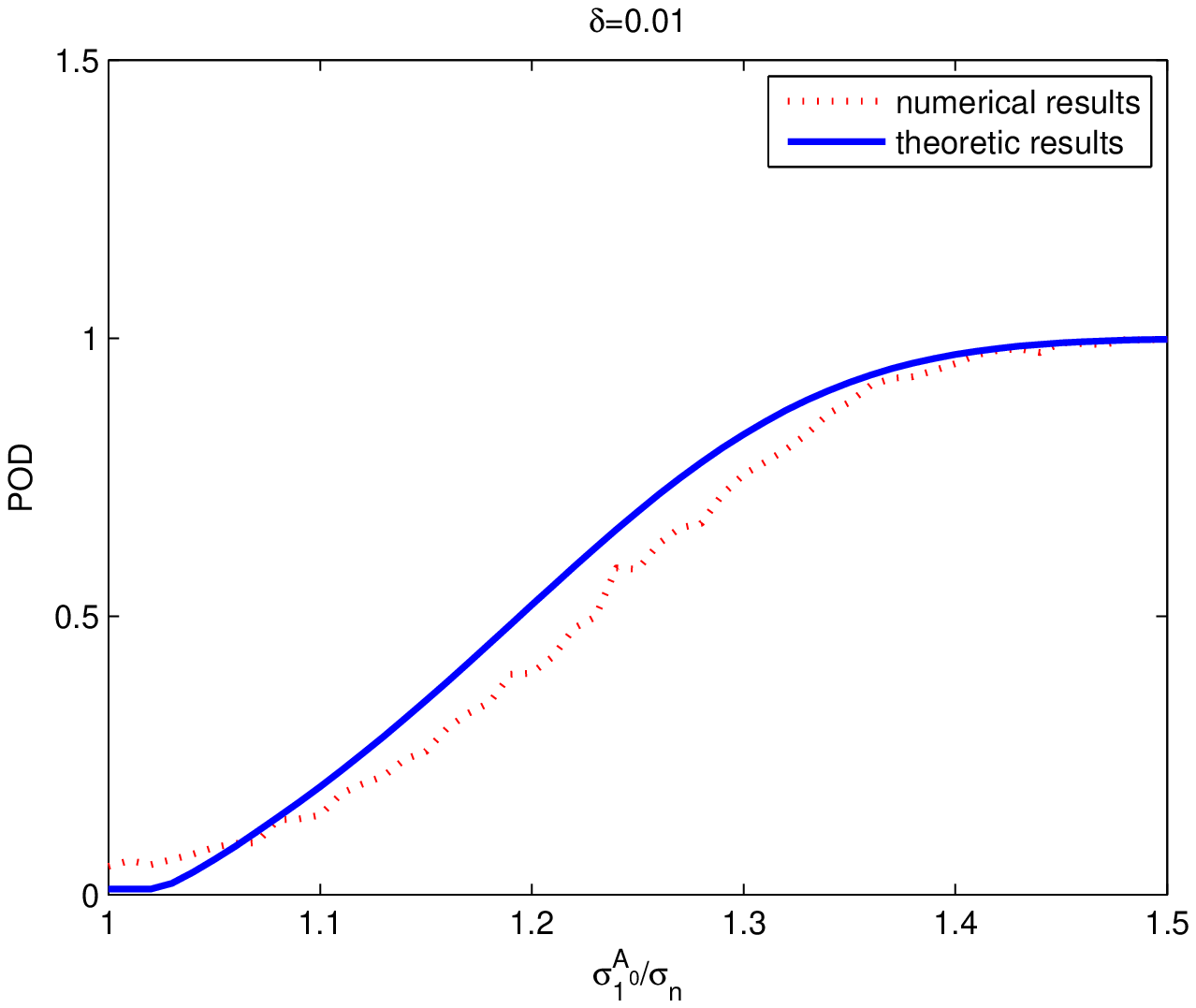}
\includegraphics[width=0.32\textwidth]{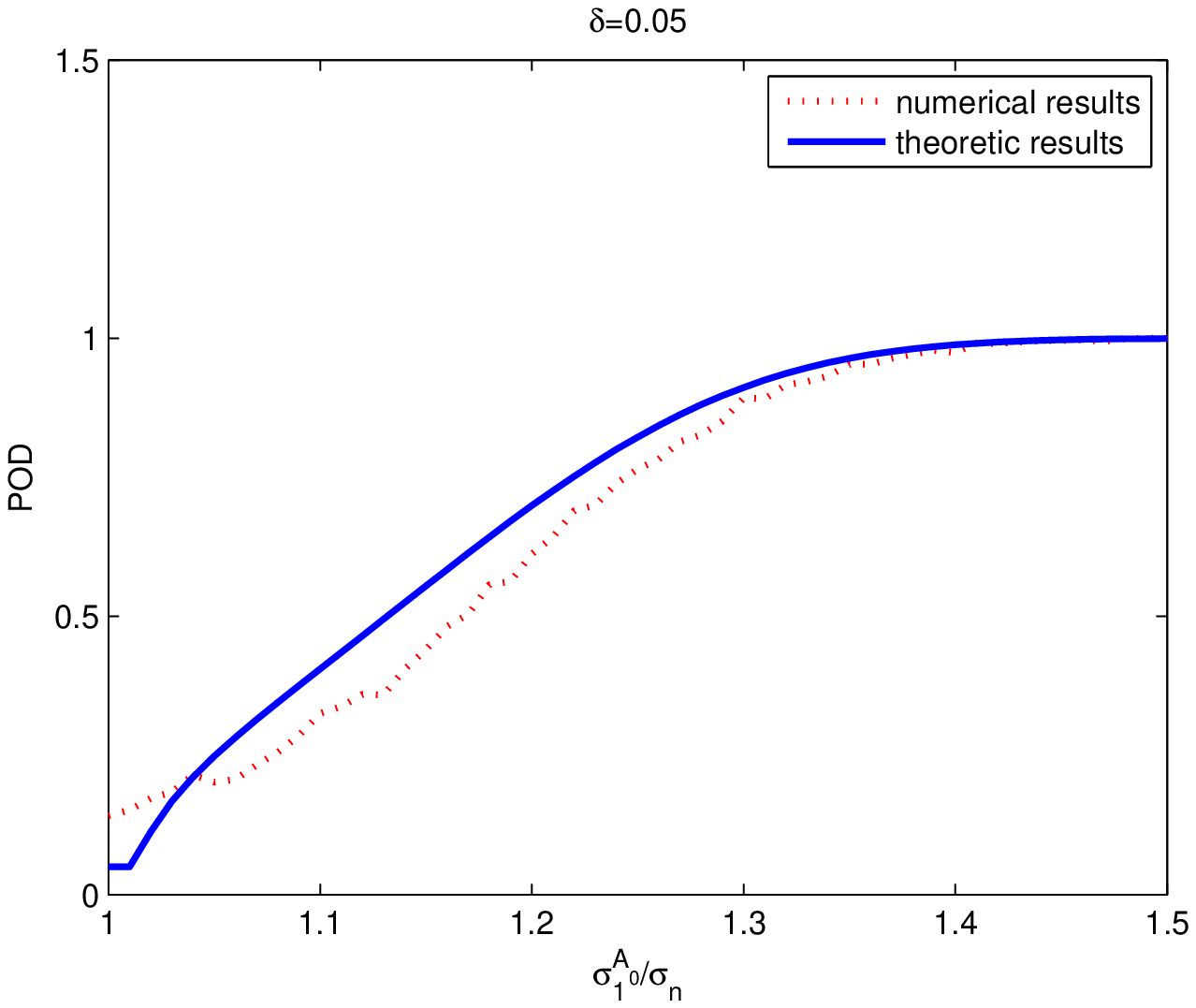}
\includegraphics[width=0.32\textwidth]{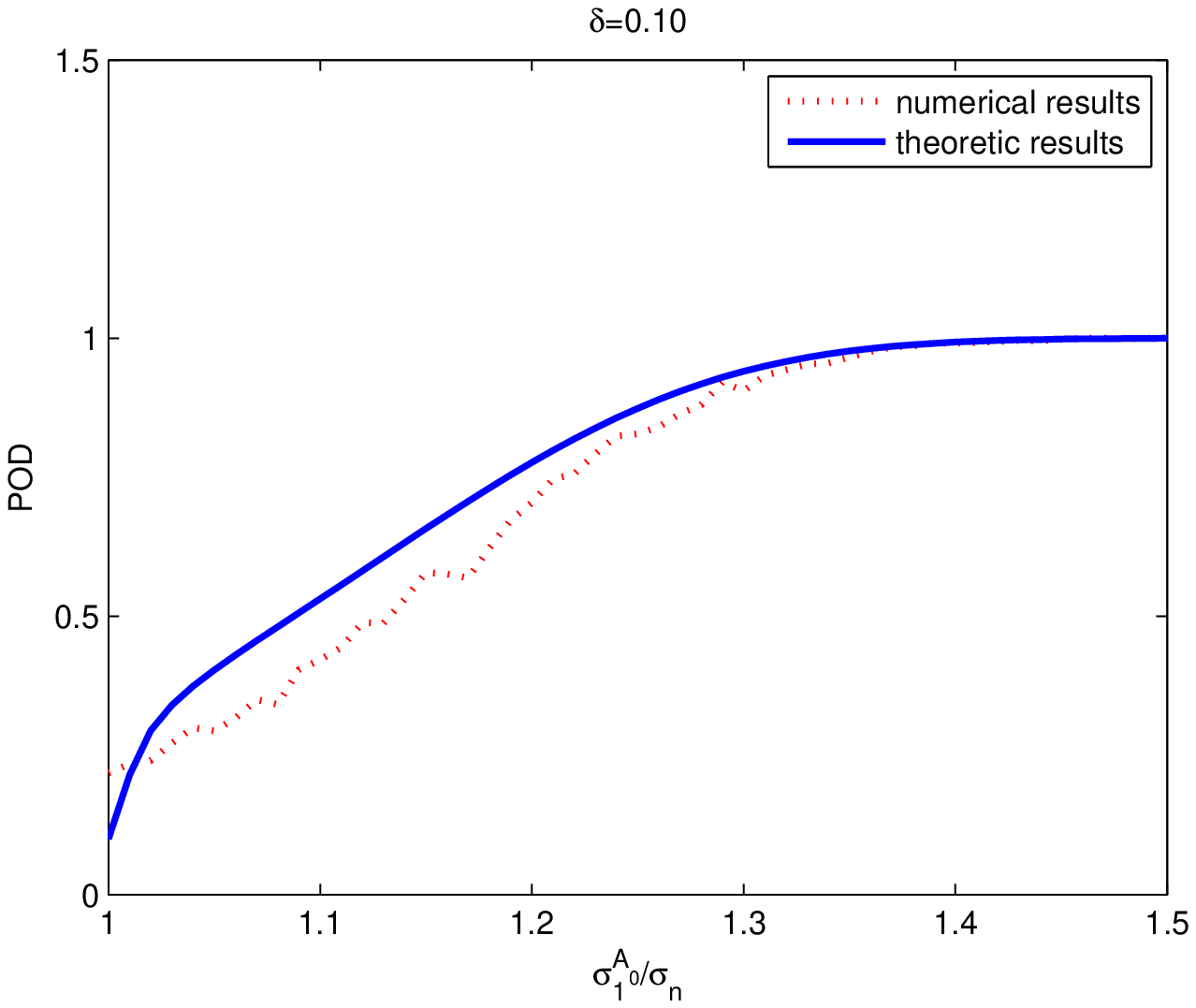}
\caption{{\rm POD} with respect to $\sigma^{\bf A_0}_1/\sigma_{\rm n}$ for different $\delta$, $\delta=0.01,
0.05, 0.10$ from left to right. }\label{fig6-3}
\end{figure}

We test the influence of the noisy measurements by adding a Gaussian noisy matrix with mean zero and variance $\sigma_{\rm n}^2/M$
to unperturbed response matrix ${\bf A}_0$. In our tests, the Gaussian noise is generated by MATLAB function {\it{randn}}. The imaging results shown in Figure \ref{fig6-2} indicate that the imaging results become sharper as the noise level is smaller.
Then we show the validity of \eqref{pod}. Noticing that $M=N$ makes $\gamma=1$ in our setting. By the analysis in Section 5, for given FAR $\delta$, POD depends on the ratio ${\sigma^{{\bf A}_0}_1}/{\sigma_{\rm n}}$. Here we only consider the critical regime in which $\sigma^{\bf A_0}_1$ is of the same order of $\sigma_{\rm n}$ (specially $\sigma^{\bf A_0}_1>\sigma_{\rm n}$).  Fixing FAR $\delta$, for each ratio ${\sigma^{\bf A_0}_1}/{\sigma_{\rm n}}$, we generate 1000  Gaussian noisy matrices with mean zero and variance $\sigma_{\rm n}^2/M$ and add them to $\bf A_0$ to get according noisy response matrices $\bf A$. We compute $R$ with the help of SVD for each $\bf A$ and count the times for $R>r_\delta$ to get the numerical POD. Figure \ref{fig6-3} shows  the comparisons between numerical POD and \eqref{pod} for each $\delta$.  We can conclude that the numerical results are in good agreement with \eqref{pod}.

\section{Concluding Remarks}

In this paper we have provided an asymptotic expansion for the
perturbations of the magnetic field due to the presence of an
arbitrary shaped small conductive inclusion with smooth boundary
and constant permeability and conductivity parameters. This was
done under the assumption that the characteristic size of the
inclusion is of the same order of magnitude as the skin depth. Our
analysis can be extended to the case of variable permeability and
conductivity distributions. We expect, however, that dealing with
nonsmooth inclusions is challenging.

Our asymptotic formula was in turn used to construct a
 method for localizing conductive targets.
  We also presented numerical simulations for
  illustration.
Thinking ahead, it appears that it would be very interesting to
apply the  findings
 from this paper to real-time
target identification in eddy current imaging using the so called
dictionary matching method \cite{dic1,dic2}.  We are also
interested in investigating target tracking from induction data at
multiple frequencies. In the presence of noise, another problem of
interest is  to study how to estimate resolution for the
localization of targets.
 This
will be the subject of a forthcoming publication.

\end{document}